\documentclass[11pt]{amsart}
\usepackage{amsmath,amsthm, amscd, amssymb, amsfonts, mathrsfs}
\usepackage[all]{xy}
\usepackage{color}
\usepackage{enumitem}

\newcommand{\nsbgp}{\vartriangleleft}
\newcommand{\Nx}{N^{[x]}}

\newcommand{\St}{\mathtt{S}}
\newcommand{\Rt}{\mathtt{R}}

\usepackage{hhline}
\usepackage[active]{srcltx}

\newcommand{\PSL}{\mathbf{PSL}}
\newcommand{\Gb}{\boldsymbol{G}}
\newcommand{\PSp}{\mathbf{PSp}}

\newcommand{\PSU}{\mathbf{PSU}}

\newcommand{\PGU}{\mathbf{PGU}}

\newcommand{\SU}{\mathbf{SU}}

\newcommand{\x}{\mathtt{x}}
\newcommand{\y}{\mathtt{y}}

\newcommand{\tJ}{\mathtt{J}}

\newcommand{\erre}{\mathtt{r}}
\newcommand{\esse}{\mathtt{s}}

\newcommand{\diag}{\operatorname{diag}}

\newcommand{\Irr}{\operatorname{Irr}}

\newcommand{\ord}{\operatorname{ord}}

\newcommand{\supp}{\operatorname{supp}}

\newcommand{\Inn}{\operatorname{Inn}}

\newcommand{\Tb}{\mathbf T}

\newcommand{\kc}{\mathbb F_q}

\newcommand\toba{\mathfrak B }

\newcommand{\trid}{\triangleright}

\newcommand{\Lb}{{\mathbb L}}
\newcommand{\Pb}{{\mathbb P}}
\newcommand{\Fc}{{\mathcal F}}

\newcommand{\qb}{{\mathbb q}}

\newcommand{\kk}{\Bbbk}

\newcommand{\ku}{\mathbb C}

\newcommand{\Z}{{\mathbb Z}}
\newcommand{\N}{{\mathbb N}}
\newcommand{\I}{{\mathbb I}}

\newcommand{\G}{{\mathbb G}}

\newcommand{\Q}{{\mathsf Q}}
\newcommand{\F}{{\mathbb F}}
\newcommand{\C}{{\mathcal C}}

\newcommand{\GL}{\mathbf{GL}}
\newcommand{\PGL}{\mathbf{PGL}}
\newcommand{\SL}{\mathbf{SL}}
\newcommand{\Sp}{\mathbf{Sp}}

\newcommand{\Oc}{{\mathcal O}}
\newcommand{\oc}{{\mathcal O}}

\newcommand{\ydg}{{}^{\ku G}_{\ku G}\mathcal{YD}}
\newcommand{\ydk}{{}^{\ku K}_{\ku K}\mathcal{YD}}

\newcommand{\Aut}{\operatorname{Aut}}

\newcommand\Ad{\operatorname{Ad}}

\newcommand{\rw}{{\widetilde r}}
\newcommand{\sw}{{\widetilde s}}

\newcommand{\clas}[2]{\mathcal C(#1,#2)}

\numberwithin{equation}{section}
\theoremstyle{plain}

\newtheorem{lema}{Lemma}[section]
\newtheorem{theorem}[lema]{Theorem}
\newtheorem{cor}[lema]{Corollary}

\newtheorem{prop}[lema]{Proposition}

\newtheorem{question-app}{Question}

\theoremstyle{definition}
\newtheorem{definition}[lema]{Definition}

\theoremstyle{remark}
\newtheorem{obs}[lema]{Remark}

\newtheorem{stepsldos}{Case}

\newcommand\id{\operatorname{id}}

\newcommand\st{\mathbb S_3}
\newcommand\sk{\mathbb S_4}

\newcommand\aco{\mathbb A_5}

\newcommand\Sim{\mathbb S}

\def\pf{\begin{proof}}
\def\epf{\end{proof}}

\theoremstyle{remark}
\newtheorem{punto}{}

\begin{document}

\renewcommand{\baselinestretch}{1.2}

\thispagestyle{empty}
\title[Finite-dimensional pointed Hopf algebras over $\PSL_n(q)$]
{Finite-dimensional pointed Hopf algebras\newline over finite simple 
groups of Lie type \newline III. Semisimple classes in $\PSL_n(q)$ }

\author[N. Andruskiewitsch, G. Carnovale, G. A. Garc\'ia]
{Nicol\'as Andruskiewitsch, Giovanna Carnovale and \\Gast\'on Andr\'es Garc\'ia}

\thanks{
This work was partially supported by
ANPCyT-Foncyt, CONICET, Secyt (UNC), the bilateral agreement between 
the Universities of C\'ordoba and Padova, the Visiting Scientist Program of the University
of Padova, the Grant CPDA 1258 18/12 of the University of Padova and an AMIDILA-ERASMUS Academic staff Fellowship.}

\address{\noindent N. A. : Facultad de Matem\'atica, Astronom\'{\i}a y F\'{\i}sica,
Universidad Nacional de C\'ordoba. CIEM -- CONICET. 
Medina Allende s/n (5000) Ciudad Universitaria, C\'ordoba,
Argentina}
\address{\noindent G. C. :
Dipartimento di Matematica,
Universit\`a degli Studi di Padova,
via Trieste 63, 35121 Padova, Italia}
\address{\noindent G. A. G. : Departamento de Matem\'atica, Facultad de Ciencias Exactas,
Universidad Nacional de La Plata. CONICET. Casilla de Correo 172, (1900)
La Plata, Argentina.}
\email{andrus@famaf.unc.edu.ar}
\email{carnoval@math.unipd.it}
\email{ggarcia@mate.unlp.edu.ar}


\begin{abstract}
We show that Nichols algebras of most simple Yetter-Drin\-feld modules over the projective 
special linear group over a finite field, 
corresponding to semisimple orbits, have infinite dimension. We introduce a new 
criterium to determine when a conjugacy class collapses and prove that for 
an infinite family of pairs $(n,q)$, any
finite-dimensional pointed Hopf algebra $H$ with  $G(H)\simeq\PSL_{n}(q)$ or $\SL_{n}(q)$
is isomorphic to a group algebra. 
\end{abstract}

\maketitle

\begin{quote}{\textit{That is not dead which can eternal lie. \\
And with strange aeons even death may die.}
}\end{quote}

\rightline{Abdul Alhazred}

\setcounter{tocdepth}{1}


\section{Introduction}

\subsection{}
This is the third paper of a series devoted to finite-dimensional
pointed Hopf algebras over $\ku$ with group of group-likes isomorphic to a finite simple group of Lie type.
An Introduction to the whole series is in Part I \cite{ACG-I}. 
Let $p$ be a prime number, $m\in \N$, $q =p^m$ and $\kc$ the field with $q$ elements.
In this paper  we consider Nichols algebras associated to semisimple conjugacy classes in $\PSL_n(q)$;
we first show that any semisimple  class $\Oc$ lying in a large family \emph{collapses} \cite[Definition 2.2]{AFGV-ampa},
that is, the dimension of the Nichols algebra $\toba(\Oc, \qb)$ is infinite for
every finite faithful 2-cocycle $\qb$. In previous work 
\cite{ACG-I, ACG-II, AFGV-ampa, AFGV-espo} we attacked the question of the collapse
of conjugacy classes in various groups using the criteria of type D and F, which are 
based on results on Nichols algebras of Yetter-Drinfeld modules over finite groups \cite{AHS, CH, HS}.
Here the criterium of type C, conjectured some years ago by the authors of \cite{AFGV-ampa, HS} 
but carried into effect only now thanks to the remarkable classification in \cite{HV}, is added to the panoply. 
See Theorem \ref{th:typeC}. Note that there are racks of type C which are not of type D or F, see 
for example Lemma \ref{lem:sl32-21-C}.

If $n=2$, we assume that $q\neq 2$, $3$, $4$, $5$, $9$
to avoid coincidences with cases treated elsewhere,
see \cite[Subsection 1.2]{ACG-I}. Indeed
$\PSL_{2}(2) \simeq \mathbb{S}_{3}$, $\PSL_{2}(3) \simeq \mathbb{A}_{4}$ are not simple and
$\PSL_{2}(4) \simeq \PSL_{2}(5) \simeq \mathbb{A}_{5}$
and  $\PSL_{2}(9) \simeq \mathbb{A}_{6}$.
Our first main result says

\begin{theorem}\label{thm:psl2-collapse}
Let $\Oc$ be a semisimple conjugacy class in $\PSL_n(q)$. If either $n>2$ and $\Oc$ is not 
irreducible or $n=2$, $q\neq 2,3,4,5,9$ and $\Oc$ is not listed in 
Table \ref{tab:ss-psl2}, then it collapses.
\end{theorem}

\begin{table}[ht]
\caption{Kthulhu semisimple classes in $\PSL_2(q)$.}\label{tab:ss-psl2}
\begin{center}
\begin{tabular}{|c|c|c|}
\hline $q$  & class &  Remark  \\
\hline  \hline
even and & irreducible & sober\\
not a square & order 3 & \\
\hline 
all & irreducible, order $> 3$ & sober  \\
\hline
\end{tabular}
\end{center}
\end{table}

\pf
Let $\Oc$ be a semisimple conjugacy class in $\PSL_n(q)$.
If $n>2$ and $\Oc$ is not
irreducible, then Proposition \ref{prop:reduction} applies.
If $n=2$, then the result follows by Proposition \ref{prop:sl2-ss}. 
\epf

\subsection{} In the first two papers of the series, we dealt with 
unipotent classes in $\PSL_n(q)$ and $\PSp_n(q)$. 
The outcome is that most non-semisimple conjugacy classes collapse, 
and yet unpublished results on other finite simple groups of Lie type
convey to the idea that this is the case in general. On the contrary, 
we see in the present paper that a semisimple irreducible conjugacy  
class might not satisfy the criteria of types C, D or F for a rack to collapse;
it appears to us that this 
would be true in general. An intuitive explanation might be as follows. 
If $G$ is a finite (almost) simple group, then there exists a  conjugacy 
class $\C$ of G so that if $x\in G - \{e\}$, 
then the probability that $x$ and a random element of $\C$ generate 
$G$ is at least 1/10; in particular $G$ can be generated
by a pair of elements in $\C$. See  \cite{BGK, GK, GM} and references therein; 
for $G$ of Lie type, $\C$ is semisimple.
Of course this does not prove that generically the conjugacy 
classes do not collapse-- see \cite{AFGV-ampa, AFGV-espo} 
for alternating and sporadic groups-- but it might be an indication 
of the plausibility of our guess.

\subsection{}
The fact that large families of conjugacy classes of finite 
simple groups do not collapse shows the 
limits of the criteria of types D, C or F, and urges for the computation 
of their second cohomology groups and the determination of the
corresponding Nichols algebras. For cocycles coming from Yetter-Drinfeld 
modules over $\PSL_n(q)$, abelian techniques can be applied.
This type of techniques were applied in \cite{FGV1} and \cite{FGV2} for the study of 
Nichols algebras over $\SL_{2}(q)$ and $\PSL_2(q)$, respectively, where it was proved that 
$\dim \toba(\Oc, \rho) = \infty$ for any irreducible representation $\rho$ of
$C_{\PSL_2(q)}(x)$ with $x\in \Oc$ and $q$ even. For $q$ odd, 
a list of the open cases was given. With the use of the criterium of type C, we deal with these open cases 
when $q\equiv 1 \mod 4$,  see Theorem \ref{thm:sl2}. 
For general $n> 1$, the use of the criteria and the abelian
techniques, together with  \cite[Theorem 1.4]{ACG-I},
yield the following result.

Let $n=2^ab$, $q=1+2^c d$ with $(b,2)=(d,2)=1$ and 
let $\mathcal{F}$ be the set of all pairs $(n,q)$ with  $n\in \N$ and $q=p^{m}$, $p$ prime, such that 
one of the following conditions holds:
\begin{enumerate}[leftmargin=*,label=\rm{(\alph*)}]
\item $n>3$ is odd;
\item $n>3$ and $q$ is even;
\item $n=3$ and $q>2$;
\item $0<a<c$, $n>2$; 
\item $a=c>1$.
\end{enumerate}

We point out that, for every fixed $n>3$, there can be only finitely many $q$ such that $(n,q)\not\in \mathcal{F}$.

\begin{theorem}\label{thm:YD-modules}Let $\Gb=\PSL_n(q)$.\begin{enumerate}
\item If $(n,q) \in \mathcal{F}$, then 
$\Gb$ collapses. Under these assumptions, any
finite-dimensional pointed Hopf algebra $H$ with  $G(H)\simeq\Gb$
is isomorphic to $\ku \Gb$. 
\item Let $n=2$, $q>3$. If  
$\Oc_x$ is not a semisimple irreducible conjugacy class which is the image of 
 $\Oc_{\x}^{\SL_{2}(q)}$ with $\x = \left(\begin{smallmatrix}
a& \zeta b\\
b & a
\end{smallmatrix}\right)$ for some $a,\,b\neq 0$, $\zeta \in \F_{q}^{\times} - \F_{q}^{2}$ and  
$q\equiv 3 \mod 4$,  then $\dim\toba(\Oc,\rho)=\infty$ for every
$\rho\in \Irr C_{\Gb}(x)$. 
\end{enumerate}
\end{theorem}

\pf
If $n=2$ we may assume that $q>5$ and $q\neq 9$, for if
$q=4,5$ or $9$, then the result follows by \cite[Theorem 1.2]{AFGV-ampa}.
Let $\Oc$ be a conjugacy class in $\Gb$ and $x\in \Oc$.
If $x$ is not semisimple, then by Proposition \ref{prop:sln-nonss}  it collapses unless $\Oc$ is
unipotent and is listed in Table \ref{tab:unip-psln2}. If it is unipotent and  
$n=2$, then $\dim \toba(\Oc, \rho) = \infty$ for all $\rho \in \Irr_{C_{\Gb}(x)}$ by \cite[Theorem 1.6]{FGV2}. 
Thus, the only open case is when $n=3$, $q=2$ and $\Oc$ is regular.
Assume now that $\Oc$ is semisimple. If $n=2$, then the result follows from Theorem \ref{thm:sl2}. Suppose
then that $n>2$. If $\Oc$ is not irreducible, then it collapses by Theorem \ref{thm:psl2-collapse}. 
Finally, if $\Oc$ is irreducible, then the result follows from Lemma \ref{lem:irred-psll-typeB}.
\epf

\subsection*{Notation} We denote the cardinal of a set $X$ by $\vert X\vert$. 
If $k < \ell$ are positive integers, then we set $\I_{k, \ell} = \{i\in \N: k\le i\leq \ell\}$ and simply
$\I_{\ell} = \I_{1, \ell}$. We denote by ${\mathbb{G}'}_\ell$ the set 
of non-trivial $\ell$-th roots of unity in ${\mathbb C}$. 

Let $G$ be a group; $N < G$, respectively $N \nsbgp G$, means that 
$N$ is a subgroup, respectively a normal subgroup, of $G$. 
The set of isomorphism classes of irreducible representations of $G$ is denoted 
by $\Irr G$. If $G$ 
acts on a set $X$ and $x\in X$ we denote by
$\Oc_x^G$ the orbit of $x$ for this action. In particular, if $x\in G$, then $\Oc_x^G$ 
indicates the conjugacy class of $x$ in $G$. 
The centralizer, respectively the normalizer, of $x\in G$ is denoted 
by $C_{G}(x)$, respectively $N_{G}(x)$;
the inner automorphism defined by conjugation by $x$ is denoted 
by $\Ad x$. If $\Fc\in \Aut G$, then $G^{\Fc}$
denotes the subgroup consisting of points fixed by $\Fc$. 
The standard Frobenius morphism will be indicated by $F_q$. 

The algebraic closure of $\F_q$ is denoted $\kk = \overline{\F_q}$. For $a,b\in{\mathbb N}$ we will set 
$(a)_b=b^{a-1}+\cdots+1$. We recall that for $a,b,c\in{\mathbb N}$ there holds
\begin{equation}\label{eq:q-number1}((a)_b, b-1)=(a, b-1),\quad((a)_b,(c)_b)=((a,c))_b,\quad
(ac)_b=(a)_b(c)_{b^a}.
\end{equation}

Let $G$ be a group and $V$ a Yetter-Drinfeld module over $G$ with comodule map $\delta$.
We denote by $V_{g} = \{v\in V:\ \delta(v)=g\otimes v\}$ the set 
of $g$-homogeneous elements for $g\in G$
and by $\supp V=\{g \in G:\ V_{g}\neq 0\}$ the support of $V$. 
Recall that  the category $\ydg$ of 
Yetter-Drinfeld modules over $G$ is braided, with 
braiding $c_{V,W}(v\otimes w)=g\cdot w\otimes v$ for $v\in V_g$,
$w\in W$, and $V,W\in \ydg$. 

\section{Preliminaries on racks}\label{sec:racks-preliminaries}

\subsection{Racks}\label{subsec:racks}

A rack is a non-empty set $X$  with a self-distributive
operation $\trid: X \times X \to X$ such that $\varphi_x := x\trid \underline{\quad}$ is a
bijection for every $x \in X$. We assume that all racks in 
this paper are finite, unless explicitly stated
and also that are crossed sets, namely that
\begin{align*}
x \trid x &= x, & x \trid y &= y \implies y \trid x = x, & \forall x,y &\in X.
\end{align*}

The main example of a rack is a conjugacy class $\Oc$
in a finite group $G$
with the operation $x\trid y = xyx^{-1}$, $x, y \in \Oc$. We say that a rack $X$ 
is \emph{abelian} if $x\trid y = y$, for all $x, y \in X$; 
thus any subset of an abelian group is an abelian rack. 

If $X$ is a rack, the inner group of the rack is ${\rm Inn} X:=\langle \varphi_x, x\in X\rangle< {\mathbb S}_X$. If 
$X=\Oc$ a  conjugacy class, then ${\rm Inn} X=\langle {\rm Ad}(y), y\in \Oc\rangle$. 

\medbreak A rack $X$ is \emph{of type D}
if it has  a decomposable subrack
$Y = R\coprod S$ with elements $r\in R$, $s\in S$ such that
$r\trid(s\trid(r\trid s)) \neq s$ \cite[Definition 3.5]{AFGV-ampa}. 
If $\Oc$ is a finite conjugacy class in $G$, then the
following are equivalent:

\begin{enumerate}
  \item\label{item:typeD-rack} The rack $\Oc$ is of type D.

  \item\label{item:typeD-group} There exist $r$, $s\in \Oc$
  such that $\Oc_r^{\langle r, s\rangle}\neq \Oc_s^{\langle r, s\rangle}$ and $(rs)^2\neq(sr)^2$.
\end{enumerate}

\medbreak A rack $X$ is \emph{of type F}
if it has a family of subracks $(R_a)_{a \in \I_4}$
and elements $r_a\in R_a$, $a \in \I_4$, such that $R_a \triangleright R_b = R_b$,  for $a, b \in \I_4$, and
$R_a \cap R_b = \emptyset$, 
$r_a\triangleright r_b \neq r_b$ for $a\neq b \in \I_4$ \cite[Definition 2.4]{ACG-I}. 
If $\Oc$ is a finite conjugacy class in $G$, then the
following are equivalent:

\begin{enumerate}
  \item\label{item:typeF-rack} The rack $\Oc$ is of type F.

  \item\label{item:typeF-group} There exist $r_a\in \Oc$, $a\in \I_4$,
  such that $\Oc_{r_a}^{\langle r_a: a\in \I_4\rangle}\neq \Oc_{r_b}^{\langle r_a: a\in \I_4\rangle}$  and
$r_a r_b \neq r_br_a$, $a  \neq b\in \I_4$.
\end{enumerate}

A rack $X$ of type D, respectively F, collapses  \cite[Theorem  3.6]{AFGV-ampa}, 
respectively \cite[Theorem 2.8]{ACG-I}.
A rack is \emph{cthulhu}, respectively \emph{sober}, 
when it is neither of type D nor of type F, respectively if  every subrack 
is either abelian or indecomposable.
Clearly, sober implies cthulhu. Let $\pi:X \to Y$ be a surjective morphism 
of racks. If $Y$ is of type D, respectively F,
then so is $X$; hence $X$ cthulhu implies $Y$ cthulhu. 

The following result extends the isogeny argument  \cite[Lemma 1.2]{ACG-I}.

\begin{lema}\label{obs:cthulhu-projection}
Let $G$ be a  group, $x\in G$, $N\nsbgp G$, $\Gb = G/N$ and 
$\pi: G \to \Gb$ the natural projection. Then the restriction
$\pi_{\oc}: \oc_x^G \to \oc_{\pi(x)}^{\Gb}$ is surjective. Assume that $G$ is finite. 
If  $\oc^G_x$ is ctuhlhu, then so is $\oc_{\pi(x)}^{\Gb}$. 
Assume that $N < Z(G)$ and let $\Nx :=\{c\in N: cx \in \oc_x^G\}$. 
Then $\Nx < N$ and $\pi_{\oc}^{-1} (\pi(y))$ has exactly
$|\Nx|$ elements. Thus, $\pi_{\oc}$ is injective if and only if $\Nx$ is trivial.
\end{lema}

\pf If $\pi(y) \in \oc_{\pi(x)}^{\Gb}$, then there is $g \in G$ 
such that $\pi(y) = \pi(g)\pi(x)\pi(g)^{-1}$ 
$= \pi(gxg^{-1}) \in \pi (\oc_x^G)$. That is, $\pi_{\oc}$ is surjective. 
It is easy to see that $\Nx$ is a subgroup
of $N$, provided that $N < Z(G)$. Let $g \in G$. If $h\in G$ 
satisfies $\pi(hxh^{-1}) =\pi(gxg^{-1})$, then there is
$c\in N$ such that $hxh^{-1} = c gxg^{-1} = gcxg^{-1}$, 
hence $cx = g^{-1}h xh^{-1}g \in \oc_x^G$. Conversely,
if $c\in \Nx$, say $cx = u x u^{-1} \in \oc_x^G$ for some $u\in G$, 
then  $(gu)x(gu)^{-1} = c gxg^{-1}$.
\epf

\subsection{Racks of type C}\label{subsec:typeC}

We now translate the main result of \cite{HV} to the context of racks, yielding a new criterium. First we recall:

\begin{theorem}\label{th:HS} \cite[Theorem 2.1]{HV}
    Let $G$ be a non-abelian group and $V$ and $W$ be two  simple Yetter-Drinfeld
    modules over $G$ such that $G$ is generated by the support of $V\oplus W$, $\dim V \le \dim W$ and 
    $(\id-c_{W,V}c_{V,W})(V\otimes W) \neq 0$. Then the following are equivalent:
    \begin{enumerate}[leftmargin=*,label=\rm{(\alph*)}]
     \item $\dim \toba (V\oplus W) < \infty$.
     \item $G$, $V$ and $W$ are as in \cite[Theorem 2.1]{HV}.
    \end{enumerate}
In particular,  $(\dim V, \dim W)$ belongs to
\begin{align}\label{eq:rank2-set}
\{(1,3), (1,4), (2,2), (2,3), (2,4)\}. \qed
\end{align}
\end{theorem}

The theorem above motivates the following definition.
\begin{definition}\label{def:rack-typeC} A rack
$X$ is \emph{of type C}  when there are  a decomposable subrack
$Y = R\coprod S$  and elements $r\in R$, $s\in S$ such that
\begin{align}
\label{eq:typeC-rack-inequality}
&r \triangleright s \neq s,
\\ \label{eq:typeC-rack-indecomposable}
&R = \Oc^{\Inn Y}_{r}, \qquad S = \Oc^{\Inn Y}_{s},
\\
\label{eq:typeC-rack-dimension} &\min \{\vert R \vert, \vert S \vert \}  > 2 \quad \text{ or } \quad
\max \{\vert R \vert, \vert S \vert \}  > 4.
\end{align}
\end{definition}

\begin{obs} Since $X$ is a crossed set, \eqref{eq:typeC-rack-inequality} 
implies that $s \triangleright r \neq r$, hence $\vert R \vert \neq 1$
and $\vert S \vert \neq 1$. That is, \eqref{eq:typeC-rack-dimension} 
says that either $\vert R \vert \neq 2$ or
 $\vert S \vert >4$.

Assume that $R$ is indecomposable and $Y=R\coprod S$ a decomposable rack. Then 
$R = \oc_r^{\Inn R} = \oc_r^{\Inn Y}$ by \cite[Lemma 1.15]{AG-adv}. 
The formulation \eqref{eq:typeC-rack-indecomposable} is
more flexible, see Lemma \ref{lema:typeC-projection}. 
On the other hand, racks with 2 elements are not indecomposable.
Thus, in presence of \eqref{eq:typeC-rack-inequality}, 
the hypothesis \emph{$R$ and $S$ are indecomposable} 
implies both  \eqref{eq:typeC-rack-indecomposable}
and \eqref{eq:typeC-rack-dimension}.
\end{obs}

\begin{obs}\label{obs:typeC-decomposition} Let $H$ be a finite group, $r\in H$ and $s\in \Oc_r^H$.
If $R:=\Oc_r^{\langle r,s \rangle}\neq S:=\Oc_s^{\langle r,s \rangle}$, 
then  $Y:=R\coprod S$ is a decomposable subrack of $\Oc_r^H$
that satisfies \eqref{eq:typeC-rack-indecomposable}, because  $\langle Y\rangle = \langle r,s\rangle$ so 
$R=\Oc_r^{\langle Y\rangle}=\Oc_r^{\Inn Y}$. 

Conversely, any subrack decomposition $Y=R \coprod S\subset \Oc_r^H$ with $r\in R$ and $s\in S$ implies
$\Oc_r^{\langle r,s \rangle}\neq \Oc_s^{\langle r,s \rangle}$ 
and \eqref{eq:typeC-rack-indecomposable} 
is verified for $R':=\Oc_r^{\langle r,s \rangle}$, $S':=\Oc_s^{\langle r,s \rangle}$. 
\end{obs}

\begin{obs}\label{rmk:genC}
Let $H$ be a finite group, $x_{1},\ldots, x_{n} \in H$ and $L= \langle 
x_{1},\ldots, x_{n}\rangle$. If $Y = \Oc_{x_{1}}^{L}\cup \cdots \cup \Oc_{x_{n}}^{L}$,
then $Y = \Oc_{x_{1}}^{\Inn Y}\cup \cdots \cup \Oc_{x_{n}}^{\Inn Y}$. Indeed,
$\langle Y \rangle = L$, 
and consequently $\Oc_{x_{i}}^{\Inn Y}=\Oc_{x_{i}}^{\langle {\rm Ad}(y),\, y\in Y\rangle}= 
\Oc_{x_{i}}^{L}$ for all $1\leq i \leq n$.    
\end{obs}

The following lemma will help us to deterine when a rack is of type C.

\begin{lema}\label{lema:typeC-odd} Let $H$ be a finite group, $r\in H$
and $s\in \Oc_r^H$ such that \eqref{eq:typeC-rack-inequality} holds.
\begin{enumerate}
[leftmargin=*,label=\rm{(\alph*)}]
  \item\label{item:typeC-odd-uno}
If $\vert \Oc_r^{\langle r,s \rangle}\vert =2$, then $s^{-1}\trid r= s\trid r$, i. e., $s^2 r=r s^2$.

\item\label{item:typeC-odd-dos}
If $\ord r = \ord s$ is odd, then  $\Oc_r^H$ is of type C if and only if 
 $\Oc_r^{\langle r,s \rangle}\neq \Oc_s^{\langle r,s \rangle}$.
  \end{enumerate}
\end{lema}

\pf \ref{item:typeC-odd-uno} Indeed, $ s^{-1} \trid r\neq r\neq s \trid r$ and  
$ s^{-1}\trid r, r,  s\trid r \in \Oc_r^{\langle r,s \rangle}$. Thus $\vert \Oc_r^{\langle r,s \rangle}\vert =2$
implies $s^{-1}\trid r= s\trid r$. \ref{item:typeC-odd-dos}
If $\ord r$ is odd, then \eqref{eq:typeC-rack-inequality}  forces $s^2 r\not=r s^2$ and 
$r^2 s\not=s r^2$, thus \eqref{eq:typeC-rack-dimension} holds; while  
\eqref{eq:typeC-rack-indecomposable} holds by 
Remark \ref{obs:typeC-decomposition}. 
\epf

We describe now the criterium of type C in group-theoretical terms.
\begin{lema}\label{lem:equivC}
Let $\Oc$ be a conjugacy class in a finite group $G$. $\Oc$ is of type C if and only 
if there are $H<G$, $r,s\in H\cap\Oc$ such that
\begin{align}\label{eq:equivC1}
&rs\neq sr;\\ 
\label{eq:equivC2} &\Oc_r^H\neq \Oc_s^H; \\
\label{eq:equivC3} &H=\langle \Oc_r^H,\Oc_s^H\rangle;\\
\label{eq:equivC4} & 
\min \{\vert\Oc_r^H\vert, \,\vert\Oc_s^H\vert\}>2 \quad\text{ or }\quad
\max \{\vert\Oc_r^H\vert, \,\vert\Oc_s^H\vert\}>4. 
\end{align} 
\end{lema}

\pf If $\Oc$ is of type C with decomposable subrack $Y=R\coprod S$, then we take 
$H=\langle Y\rangle$. Conversely, if $H$, $r$ and $s$ satisfy 
\eqref{eq:equivC1} and \eqref{eq:equivC2}, then we take 
$R=\Oc_r^H$, $S=\Oc_s^H$ and $Y=R\coprod S$. Thus, $\langle Y\rangle =H$ by \eqref{eq:equivC3} hence 
\eqref{eq:typeC-rack-indecomposable} is satisfied. Finally, \eqref{eq:typeC-rack-dimension} follows
from \eqref{eq:equivC4}.
\epf

The authors of \cite{AFGV-ampa} had in mind a version of the criterium of type C, 
inspired by personal communications of the authors of
\cite{HS}, who conjectured a version of Theorem \ref{th:HS}.  
At the time C stood for conjectural.

\begin{theorem}\label{th:typeC}
A rack $X$ of type C collapses.
\end{theorem}

\pf
Let $G$ be a finite group and  $M\in\ydg$ such that $X$ is isomorphic to a
subrack of $\supp M$. We will check that $\toba(M)$ has infinite dimension. 
This implies that $X$ collapses by \cite[Lemma 2.3]{AFGV-ampa}.

Let $Y = R\coprod S$ be as in Definition \ref{def:rack-typeC}.
Let $K = \langle Y\rangle \le G$.
Then $M_Y := \oplus_{y\in Y} M_y \in \ydk$, with
$M_R := \oplus_{x\in R} M_x$ and $M_S := \oplus_{z\in S} M_z$ being  
Yetter-Drinfeld submodules of $M_Y$.
By \eqref{eq:typeC-rack-indecomposable}, $R = \Oc^{K}_{r}$, $S = \Oc^{K}_{s}$. 
Let $V$, respectively $W$, be a simple Yetter-Drinfeld submodule of 
$M_R$, respectively $M_S$. Then
$\supp V = R$ (since $\supp V$ is stable under the conjugation of $K$), 
$\supp W = S$ and $\supp (V \oplus W) = Y$, that generates $K$.
Now $(\id-c_{W,V}c_{V,W})(V\otimes W) \neq 0$ by \eqref{eq:typeC-rack-inequality}.
Without loss of generality, we may assume that $\dim V \leq \dim W$.
Now $\dim V \ge \vert R \vert > 2$ or $\dim W \geq \vert S\vert >4$, by \eqref{eq:typeC-rack-dimension}.
Hence $(\dim V, \dim W)$ does not belong to the set \eqref{eq:rank2-set}.
Thus $\dim \toba(V \oplus W) = \infty$ by Theorem \ref{th:HS} and a fortiori $\dim \toba(M) = \infty$.
\epf

The following lemma will allow us to use inductive reasoning when applying the criterium.

\begin{lema}\label{lema:typeC-projection}
If a rack $Z$ contains a subrack of type C, respectively projects 
onto a rack of type C, then $Z$ is
of type C. 
\end{lema}

\pf The first statement is obvious.  Let $\pi:Z \to X$ be a surjective morphism of racks with $X$ of type C
and let $Y = R \coprod S\subset X$ be as in Definition \ref{def:rack-typeC} with 
$\vert R \vert \leq \vert S \vert$; in particular, $\vert R \vert > 2$ or $\vert S \vert > 4$. Fix $\rw, \sw \in Z$ such that
$\pi(\rw) = r$, $\pi(\sw) = s$.
Define recursively
\begin{align*}
 R_1 &= \pi^{-1} (R),& S_1 &= \pi^{-1} (S),& Y_1 &= \pi^{-1} (Y),& K_1 &= 
 \langle \varphi_y, y\in Y_1\rangle \le \Inn Z,\\ 
 R_2 &= \Oc^{K_1}_{\rw},& S_2 &= \Oc^{K_1}_{\sw},& Y_2 &= R_2 
 \coprod S_2,& K_2 &= \langle \varphi_y, y\in Y_2\rangle \le \Inn Z; \\
 R_{j} &= \Oc^{K_{j-1}}_{\rw},& S_{j} &= \Oc^{K_{j-1}}_{\sw},& Y_{j}
 &= R_{j} \coprod S_{j},& K_{j} &= \langle \varphi_y, y\in Y_j\rangle \le \Inn Z.
 \end{align*}
Notice that $R_1 \supseteq R_2 \supseteq \dots$ and 
$S_1 \supseteq S_2 \supseteq \dots$, hence $Y_{i} = R_{i} \coprod S_{i}$
is a rack decomposition. Now the sequence $Y_1 \supseteq Y_2 \supseteq \dots \supseteq
Y_i \supseteq Y_{i+1} \supseteq \dots$
stabilizes because $Z$ is finite. Let $i \in \N$ such that $Y_i = Y_{i - 1}$; then 
$\widetilde R := R_i = R_{i-1} = \Oc^{K_{i-1}}_{\rw}$ and 
$\widetilde S := S_i = S_{i-1} = \Oc^{K_{i-1}}_{\sw}$. 
Thus $\widetilde Y := \widetilde R \coprod \widetilde S$ is a subrack of $Z$
that satisfies \eqref{eq:typeC-rack-inequality} and \eqref{eq:typeC-rack-indecomposable}. 
We claim now that $\pi (Y_j) = Y$ for all $j\in \N$; hence 
$\vert R_j\vert \ge \vert R\vert > 2$ or  $\vert S_j\vert \geq \vert S \vert  > 4$, proving 
\eqref{eq:typeC-rack-dimension} for $\widetilde Y$. 

Indeed, $\pi (R_1) = R$ because $\pi$ is surjective. Assume that $\pi (Y_j) = Y$; hence $\pi (R_j) = R$ 
and $\pi (S_j) = S$. Let $t \in R$. There exist $y_1, \dots, y_h\in Y$ such that 
$y_1\triangleright(y_2 \triangleright \dots \triangleright (y_h \triangleright r)\dots ) = t$
by \eqref{eq:typeC-rack-dimension} for $Y$. Pick  $\widetilde y_1, \dots, \widetilde  y_h\in Y_j$ such that
$\pi (\widetilde  y_\ell) = y_\ell$, $\ell\in \I_h$. Then 
\begin{align*}
 \widetilde y_1\triangleright(\widetilde y_2 \triangleright \dots 
 \triangleright (\widetilde y_h \triangleright \widetilde r)\dots ) 
 &\in \Oc^{K_j}_{\rw} = R_{j+1}, \qquad \text{hence}\\ 
 \pi(\widetilde y_1\triangleright(\widetilde y_2 \triangleright \dots 
 \triangleright (\widetilde y_h \triangleright \widetilde r)\dots ) 
 &= y_1\triangleright(y_2 \triangleright \dots \triangleright (y_h \triangleright r)\dots ) = t \in \pi(R_{j+1}).
 \end{align*}
\epf

\subsection{Kthulhu racks}\label{subsec:kthulhu}
\begin{definition}\label{def:rack-kthulhu} A rack is \emph{kthulhu} 
if  it is neither of type D nor of type F, nor of type C; i. e. 
cthulhu and not of type C. A rack is \emph{austere} if every subrack generated by two elements is 
either abelian or  indecomposable.
Clearly, sober implies austere and austere implies kthulhu. 
\end{definition}

Let $\pi:X \to Y$ be a surjective morphism of racks. 
By Lemma \ref{lema:typeC-projection} and previous results,
$X$ kthulhu implies $Y$ kthulhu.

In \cite{ACG-I} and \cite{ACG-II} we proved that the non-semisimple classes in $\PSL_n(q)$ that are not listed in Table 
\ref{tab:unip-psln} and the unipotent classes in $\PSp_{2n}(q)$ that are not listed in Table \ref{tab:unip-simplectic}
are either of type D or F. In this Section we determine which ones are of type C. 

\begin{table}[ht]
\caption{Unipotent classes in $\PSL_{n}(q)$ not of type D.}\label{tab:unip-psln}
\begin{center}
\begin{tabular}{|c|c|c|c|c|}
\hline $n$  & type & $q$ & Remark & kthulhu \\
\hline
$2$  & $(2)$ & even or & sober&  yes\\
& & not a square  & \cite[Lemma 3.5]{ACG-I} & \\
\hline
$3$ &$(3)$  & 2  & sober, \cite[Lemma 3.7 (b)]{ACG-I}& yes  \\
\cline{2-5}
& $(2,1)$ & 2 &  cthulhu, \cite[Lemma 3.7 (a)]{ACG-I}  & \\
\cline{3-4}
&  & even $\ge 4$ & cthulhu, \cite[Prop. 3.13, 3.16]{ACG-I} & type C\\
\cline{1-4}
$4$ & $(2,1,1)$  & 2  & cthulhu, \cite[Lemma 3.12]{ACG-I} & \\
\cline{3-4}
&  & even $\ge 4$ & not type D, \cite[Prop.	 3.13]{ACG-I} & Lemma \ref{lem:sl32-21-C} \\
& & & open for type F & \\
\hline
\end{tabular}
\end{center}
\end{table}

\begin{table}[ht]
\caption{Cthulhu unipotent classes in $\PSp_{2n}(q)$.}\label{tab:unip-simplectic}
\begin{center}
\begin{tabular}{|c|c|c|c|c|}
\hline $n$  & type  & $q$ & Remark & kthulhu \\
\hline
$\geq 2$  & $W(1)^a\oplus V(2)$  & even   & cthulhu  & yes
\\
\cline{2-3}
&  $(1^{r_1}, 2)$  & odd, 9 or  &\cite[Lemma 4.22]{ACG-II}  & Lemma \ref{lem:21notC}\\
&    &  not a square &  & 
\\
\hline
$3$ & $W(1)\oplus W(2)$ & 2 & cthulhu& type C \\
& &  &\cite[Lemma 4.25]{ACG-II} & Lemma \ref{lem:W(2)W(1)C} \\
\hline
$2$ & $W(2)$  & even  & cthulhu  & yes \\
&    &   & \cite[Lemma 4.26]{ACG-II} & Lemma \ref{lem:21notC}
\\
\cline{2-5}
& $(2, 2)$ & 3 &   one class cthulhu & type C\\
&  &  & \cite[Lemma 4.5]{ACG-II} &  Lemma \ref{lem:sp43-22-C}\\
\cline{2-5}
& $V(2)^2$ & 2 & cthulhu  &  type C\\
& & & \cite[Lemma 4.24]{ACG-II} & Lemma \ref{lem:22q2-1-C}\\
\hline
\end{tabular}
\end{center}
\end{table}

We recall that by the isogeny argument 
in \cite[Lemma 1.2]{ACG-I}, for unipotent classes we can work 
in $G=\SL_n(q)$ and $G=\Sp_{2n}(q)$.

\begin{lema}\label{lem:sl32-21-C}
Let $G = \SL_{n}(q)$ with $n\geq 3$ and $q$ even. 
Then any unipotent conjugacy class $\Oc$ with associated partition
$(2,1^{n-2})$ is of type C. 
\end{lema}

\pf
By Lemma \ref{lema:typeC-projection}, it is enough to prove the 
assertion for $G= \SL_{3}(2) $. 
Denote by $\alpha_{1}$ and 
$\alpha_{2}$ the positive simple roots and 
let $x=x_{\alpha_{1}+\alpha_2}(1)= \left(\begin{smallmatrix}
              1 & 0 & 1\\
	     0 & 1 & 0\\
	     0 & 0 & 1
             \end{smallmatrix}\right)
$,
$y= \left(\begin{smallmatrix}
              1 & 0 & 0\\
	     0 & 0 & 1\\
	     0 & 1 & 0
             \end{smallmatrix}\right)
$ and 
$z=x_{-\alpha_2}(1)= \left(\begin{smallmatrix}
              1 & 0 & 0\\
	     0 & 1 & 0\\
	     0 & 1 & 1
             \end{smallmatrix}\right)
$. Then $x,y$ and $z$ are conjugated in $\SL_{3}(2)$  
with $y = v\trid x$, $z= w\trid x$ and 
$v= \left(\begin{smallmatrix}
              0 & 1 & 0\\
	     1 & 0 & 1\\
	     1 & 0 & 0
             \end{smallmatrix}\right)
$, 
$w= \left(\begin{smallmatrix}
              0 & 1 & 0\\
	     0 & 0 & 1\\
	     1 & 0 & 0
             \end{smallmatrix}\right)
$.
Denote by $H$ the subgroup of 
$G$ generated by $x,y$ and $z$.
As $x, y$ and $z$ lie in the $F_q$-stable parabolic subgroup $\Pb$ of $\SL_n(\kk)$ 
with $F_q$-stable Levi factor $\Lb$ 
with root system $\{\pm \alpha_2\}$, it follows that 
$H\subset \Pb^{F_q}\subsetneq G$. Moreover, since $x$ is in the 
unipotent radical of $\Pb$, which is normal, and $y\in  \Lb^{F_q}$ 
we have that $\Oc_{x}^{H} \neq \Oc_{y}^{H}$ and a direct 
computation shows that $\Oc_{x}^{H}=\{x,
 \left(\begin{smallmatrix}
              1 & 1 & 0\\
	     0 & 1 & 0\\
	     0 & 0 & 1
             \end{smallmatrix}\right), \left(\begin{smallmatrix}
              1 & 1 & 1\\
	     0 & 1 & 0\\
	     0 & 0 & 1
             \end{smallmatrix}\right)\}$ and $\Oc_{y}^{H}=\Oc_{z}^{H}$.
Thus $|\Oc_{x}^{H}|=3$. The class $\Oc_y^H$ contains $y, z$ and $y\trid z=\ ^t\!z\neq z,y$. 
The result follows by
Remark \ref{rmk:genC} and 
Lemma \ref{lem:equivC}.
\epf

As a consequence of the lemma above, 
Theorem \ref{th:typeC} and \cite[Theorem 1.3]{ACG-I} we obtain the
following.

\begin{prop}\label{prop:sln-nonss}Let $x\in \Gb$ and pick 
$\x\in \SL_{n}(q)$ such that $\pi(\x) = x$, with Jordan
decomposition $\x = \x_s\x_u$.  Assume that $\x_u\neq e$. 
Then either  $\Oc = \oc^{\Gb}_{x}$ collapses or else $\x_s$ is central
and $\Oc$ is a unipotent class listed in Table \ref{tab:unip-psln2}. \qed
\end{prop}

\begin{table}[ht]
\caption{Kthulhu unipotent classes in $\PSL_{n}(q)$.}\label{tab:unip-psln2}
\begin{center}
\begin{tabular}{|c|c|c|c|}
\hline $n$  & type & $q$ & Remark \\
\hline
$2$  & $(2)$ & even or not a square  & sober,  \cite[Lemma 3.5]{ACG-I} \\
\hline
$3$ &$(3)$  & 2  & sober,    \cite[Lemma 3.7]{ACG-I} (b)  \\
\hline
\end{tabular}
\end{center}
\end{table}

For the rest of the Section $G=\Sp_{2n}(q)$, $\Gb=\PSp_{2n}(q)$, 
$n\geq2$ and $u\in G$ is a unipotent element.
Recall that $\clas{G}{u}$ denotes the set of $G$-conjugacy classes 
contained in $\Oc_u^{\G}$. For unexplained notation see \cite[4.2.1]{ACG-II}.

\begin{lema}\label{lem:21notC}
Let $G = \Sp_{2n}(q)$ and $\Oc$ a conjugacy class corresponding to 
\begin{enumerate}[leftmargin=*,label=\rm{(\alph*)}]
 \item  $W(1)^{a}\oplus V(2)$ or $W(2)$ if $q$ is even, or
 \item a partition $(2,1^{n-2})$ if $q$ is odd.
\end{enumerate}
Then $\Oc$ is austere, hence kthulhu.
\end{lema}

\pf  By the proof of \cite[Lemma 4.26]{ACG-II} there exists an automorphism of $\Sp_4(q)$ 
mapping each class labeled by $W(2)$ to a class corresponding to $W(1)^{2}\oplus V(2)$.
On the other hand, if $\Oc$ is labeled by $W(1)^{a}\oplus V(2)$ or 
by $(2,1^{n-2})$, then it is austere by \cite[Lemma 4.22]{ACG-II}. 
\epf

\begin{lema}\label{lem:sp43-22-C}
 Let $G = \Sp_{4}(3)$. Then the conjugacy class $\Oc$ of 
type $(2,2)$ is of type C.
\end{lema}

\pf There are two classes of type $(2,2)$, represented by $w=\left(\begin{smallmatrix}1&&&-1\\
&1&1\cr
&&1\\
&&&1\\
\end{smallmatrix}\right)$ and 
$z=\left(\begin{smallmatrix}1&&&1\\
&1&1\cr
&&1\\
&&&1\\
\end{smallmatrix}\right)$. By \cite[Lemma 4.5]{ACG-II}, $\Oc_{w}^{G}$ is of type $D$ and
$\Oc_{z}^{G}$ is chutlhu. We show that the latter is of type $C$.

Let $y=\left(\begin{smallmatrix}2&&&2\\
&1&0\cr
&1&1\\
1&&&
\\
\end{smallmatrix}\right)$, then $y \in \Oc_{z}^{G}$ since
$y = vzv^{-1}$ with $v= 
\left(\begin{smallmatrix}1&1&1&1\\
 0&0&2& 1\\
 1& 2&0&0\\
 1& 1 & 0&0 \\
\end{smallmatrix}\right)\in \Sp_{4}(3)$.
Let $H = \langle z,y \rangle \subset G$. 
Since $\ord(z)=3$, 
by Lemma \ref{lema:typeC-odd} $(b)$ it is enough to prove that 
$\oc^{H}_{z}\neq \oc^{H}_{y}$. 
Let $\mathbb{M}$ be the $F_q$-stable subgroup of $\Sp_{4}(\kk)$ of matrices 
$\left(\begin{smallmatrix}a&0&b\\
0&M&0\\
c&0&d\\
\end{smallmatrix}\right)$ with $\left(\begin{smallmatrix}a&b\\
c&d\\
\end{smallmatrix}\right)$, $M\in\SL_{2}(\kk)$. 
Clearly, $\mathbb{M}\simeq \SL_2(\kk)\times \SL_{2}(\kk)$ and  
$H \subseteq {\mathbb{M}^{F_q}} \simeq {\SL_2(3)}\times {\SL_{2}(3)}$.

Assume $y = AzA^{-1}$ with $A\in H$. Then, there
exist  $\left(\begin{smallmatrix}a&b\\
c&d\\
\end{smallmatrix}\right)$, $M\in\SL_{2}(3)$ such that 
$A=
\left(\begin{smallmatrix}a&0&b\\
0&M&0\\
c&0&d\\
\end{smallmatrix}\right)$. But then 
$\left(\begin{smallmatrix}a&b\\
c&d\\
\end{smallmatrix}\right) \left(\begin{smallmatrix}1&1\\
0&1\\
\end{smallmatrix}\right) = \left(\begin{smallmatrix}2&2\\
1&0\\
\end{smallmatrix}\right)
\left(\begin{smallmatrix}a&b\\
c&d\\
\end{smallmatrix}\right)$ and this implies that 
$a=c\neq0$, $d=2a +b$ and consequently $ad-bc= 2$, a contradiction.
Thus, $z$ and $y$ are not conjugated in $H$.\epf

Now we show that the remaining cthulhu classes in Table \ref{tab:unip-simplectic} 
are of type C.

\begin{lema}\label{lem:22q2-1-C} Assume $G=\Sp_{4}(2)$
and $\Oc$ is of the form  
$V(2)^2$. Then $\Oc$  is of type C.
\end{lema}

\pf By \cite[Theorem 6.21]{LS},
there is only one class in $\clas{G}{u}$ which is cthulhu
by \cite[Lemma 4.24]{ACG-II}.  
Since $\G^{F_q}=\Sp_4(2)\simeq\Sim_6$ and $\oc$ corresponds to the 
partition $(1^2,2^2)$, we 
show that the latter class in $\Sim_{6}$ is of type C.
Let $x=(1,2)(3,4)$, $y=(3,6)(4,5)$ and $z=(1,6)(2,5)$ and
denote by $H$ the subgroup generated by $x$, $y$ and $z$. Since 
they are all even permutations,
$H\subseteq \mathbb{A}_{6}\subsetneq \Sim_{6}$. Further, a direct computation shows that 
$\Oc_{x}^{H} = \{x, (1,2)(5,6), (3,4)(5,6)\}$, with $y\trid x =(1,2)(5,6) $, 
$z\trid x = (3,4),(5,6)$ and $z\trid(y\trid x) = x\trid(y\trid x)= y\trid x$, 
$y\trid(z\trid x) = x\trid(z\trid x)= z\trid x$. In addition,
since $(z\trid y)\trid y = z = (y\trid z)\trid y$, we have
$\Oc_{y}^{H} = \Oc_{z}^{H}$ and 
$\Oc_{y}^{H} =\{y,z,(13)(24),(35)(46),(15)(26),(14)(23)\}$. 
Hence,
$|\Oc_{x}^{H}|=3$, $|\Oc_{y}^{H}|=6$, 
$\Oc_{x}^{H}\neq \Oc_{y}^{H}$ and the result follows
by Remark \ref{rmk:genC} and Lemma \ref{lem:equivC}.
\epf

\begin{lema}\label{lem:W(2)W(1)C} 
Let $G=\Sp_{6}(2)$
and assume $\Oc$ is of the form   
$W(1)\oplus W(2)$.  
Then $\clas{G}{u}$ consists of only one class $\Oc$ which is of type C. 
\end{lema}

\pf By \cite[Lemma 4.25]{ACG-II}, $\clas{G}{u}$ consists of only one cthulhu class, 
represented by $u=x_{\alpha_1}(1)=\id_{6}+e_{1,2}+2e_{5,4}$.

Let $\tJ=\left(\begin{smallmatrix} 
                 0& 0&1\\
                 0& 1&0\\
                 1& 0&0
                 \end{smallmatrix}\right)$.
We recall that there is a natural embedding $\iota\colon \SL_3(q)\to G$ given by
$A\mapsto \diag(A, \tJ ^t\!A^{-1}\tJ)$. 
By Lemma \ref{lem:sl32-21-C}, for $x=\id_3+e_{12}$ the class $\Oc_x^{\SL_3(q)}$ is of type C.
Since $\iota(x)=x_{\alpha_1}(1)$ we have the statement. \epf

\medskip

Using the results in this section and \cite[Theorem 1.1]{ACG-II} we have the following.

\begin{prop}
Let $\Oc$ be a unipotent conjugacy class in $\PSp_{2n}(q)$.
If $\Oc$ is not listed in Table \ref{tab:unip-simplectic2}, then it collapses. \qed
        \end{prop}

\begin{table}[ht]
\caption{Kthulhu unipotent classes in $\PSp_{2n}(q)$.}\label{tab:unip-simplectic2}
\begin{center}
\begin{tabular}{|c|c|c|c|c|}
\hline $n$  & type  & $q$ & Remark \\
\hline
$\geq 2$  & $W(1)^a\oplus V(2)$  & even   & austere \\
\cline{2-3}
&  $(1^{r_1}, 2)$  & odd, 9 or  & \\
&    &  not a square & Lemma \ref{lem:21notC}
\\
\cline{1-3}
$2$ & $W(2)$  & even  & \\
\hline
\end{tabular}
\end{center}
\end{table}

\section{Preliminaries on semisimple classes in $\PSL_{n}(q)$}
\label{sec:pls-preliminaries}

From now on, $G = \SL_{n}(q)$, $\Gb = \PSL_{n}(q)$ and $\pi: G \to \Gb$ 
is the natural projection.

\subsection{Semisimple irreducible conjugacy classes}\label{subsec:irred}

Let $\St \in \GL_n(q)$ be semisimple and $\chi_{\St} = X^n + a_{n-1} X^{n-1} + \dots  
+ a_{1}  X + (-1)^{n} \in \F_q[X]$ its characteristic polynomial. 
It is well-known, see e. g. \cite[Remark 4.1]{ACG-I}, that $\oc^G_{\St} = \oc^{\GL_{n}(q)}_{\St} \cap G$.
Hence, for $\St\in\SL_n(q)$: 
\begin{equation}\label{eq:orbit-ss}
\oc^G_{\St} =\{\Rt \in  \SL_{n}(q): \Rt \text{ is semisimple and } \chi_{\Rt} =\chi_{\St}  \}.
\end{equation}

Assume that $\St  \in G$ is irreducible, that is,
$\chi_{\St}$ is irreducible. Thus 
\begin{align}\label{eq:orbit-irred}
\oc^G_{\St} =\{\Rt \in  \SL_{n}(q): \chi_{\Rt} =\chi_{\St}  \}.
\end{align}

Hence, the irreducible semisimple conjugacy classes in $G$ are 
parametrized by the monic irreducible 
polynomials of degree $n$ with constant term equal to $(-1)^n$, 
and they can be represented by the companion matrix. 

Recall \cite[Definition 3.8]{AFGV-ampa} that $g$ and its class $\Oc$ are

\begin{itemize}
  \item \emph{real} when $g^{-1}\in\oc$;
  \item \emph{quasi-real} when $g$ is not an involution and there is $j\in \Z$ such that $g\neq g^{j} \in \Oc$.
\end{itemize}


\begin{obs}\label{obs:standard} We collect some standard facts about $\oc^G_{\St}$.
\begin{punto}\label{punto:ss-a} If $\St  \in G$ is irreducible, then
the subalgebra of matrices in $M_n(q)$ commuting with $\St$  is isomorphic to $\F_{q^{n}}$. 
Hence  $C_{\SL_n(q)}(\St) \simeq \Z/(n)_q$. 
\end{punto}

\begin{punto}\label{punto:ss-b} Let $\St  \in G$ be semisimple. Then
$\chi_{\St} =\chi_{\St^q}$, hence $\St^{q} \in \Oc_{\St}^{G}$, but $\St \neq \St^{q}$ 
unless it is 
diagonalizable over $\F_q$. 
If $\eta\in\kk$ is an eigenvalue of $\St$, then $\eta^{q^l}$ is again so, for $l=1,\ldots, n-1$. 
They are all distinct if and only if $\St$ is irreducible. In this case, $\ord \St$ 
divides $(n)_q$.  If $q+1=\ord \St$, then the conjugacy class $ \Oc_{\St}^{G}$ is
\emph{real}. It is \emph{quasi-real} if it is not an involution nor diagonalizable over $\F_q$.
\end{punto}

\begin{punto}\label{obs:Cq-scalar}If $\St\in \GL_n(q)$ is semisimple 
irreducible such that 
$\St^q=\lambda \St$ for some $\lambda\in \F_q$, then $\lambda$ is a
primitive $n$-th root of $1$. In particular $n|(q-1)$. Indeed, $\St$ and $\St^q$ are conjugate, 
so $\det(\St)=\det (\St^q)=\lambda^n \det(\St)$ whence $\lambda^n=1$. In addition, 
we have $\St^{q^j}=\lambda^j \St$ for $j\in \I_{n-1}$. Since all such matrices are 
distinct, we have the claim.
\end{punto}

\begin{punto}\label{punto:ss-c}
If $\lambda\in \F_q$, $\lambda^n =1$, then 
	\begin{align*}
	 \chi_{\lambda\St} = X^n + a_{n-1}\lambda X^{n-1} + \dots + a_{j}\lambda^{n-j} 
	 X^{j} 
	 + \dots + a_{1} \lambda^{n-1} X + (-1)^{n}.
	\end{align*}
Hence, for $\St\in \GL_n(q)$  semisimple,  $\lambda\St \in \Oc_{\St}^{G}$ 
if and only if $a_{j}(1 - \lambda^{n-j}) = 0$ for every $j\in\I_{n-1}$. By \ref{obs:Cq-scalar} 
if $\St$ is irreducible, then the characteristic polynomial 
of $\St$ 
is $X^n+(-1)^n \det \St$. 
\end{punto}

\begin{punto}\label{punto:ss-d}
Assume that $q = t^h$ 
and $\St \in \SL_n(t)< \SL_n(q)$, with characteristic polynomial $\chi_{\St, t}$. 
Then $\chi_{\St, q} = \chi_{\St, t}$ (because they are determinants of the same matrix) 
and $\oc^G_{\St} \cap \SL_n(t) = \oc^{\SL_n(t)}_{\St}$ 
by \eqref{eq:orbit-ss}.
\end{punto}
\end{obs}

The following lemma will be useful in the sequel. 

\begin{lema}\label{lem:three-powers}Let $n>2$ and let 
$x\in \Gb$ be an irreducible semisimple element. 
Then $x^{q^i}\neq x^{q^j}$ for every $i\not\equiv j\mod  n$. In particular, $x$ is quasi-real. 
\end{lema} 
\pf As $x^{q^n}=x$, it is enough to prove that $x^{q^j}\neq x$ for every $j\in \I_{n-1}$.
Assume $x^{q^j}=x$. We may assume that $j\mid n$. Indeed if $x^{{q^j}-1}=1$, then 
$$\ord x\mid ((q^j-1),(n)_q) =((q-1)(j)_q,(n)_q)=(j,n)_q\left(q-1,\frac{(n)_q}{(j,n)_q}\right)$$
and the latter divides $q^{(j,n)}-1$. Let $n=jk$ and let $\x\in\SL_n(q)$ such that $\pi(\x)=x$. 
Then $\x^{q^j}=\lambda \x$ for some $\lambda\in \F_q$ and 
$\x^{q^{\ell j}}=\lambda^\ell x$ for every $\ell$. In particular, 
 $\x=\x^{q^{n}}=\lambda^k \x$ so $\lambda^k=1$. In addition, 
 as all such powers of $\x$ are distinct, $\lambda$ is a primitive $k$-th root of $1$, so $k\mid (q-1)$. 

Let $\eta$ be an eigenvalue of $\x$. Then $\eta^{q^j}=\lambda\eta$. 
Thus, the eigenvalues of $\x$ are
$\lambda^t\eta^{q^i}$ for $t\in \I_{k-1}$ and $i\in\I_{j-1}$. 
Therefore $1=\det\x=\eta^{(j)_q}\lambda^{j\binom{k}{2}}$.
Taking the $(q-1)$-th power we have
$\eta^{q^{j}-1}\lambda^{(q-1)j\binom{k}{2}}=1$. 
If $q$ is odd then $(q-1)/2\in{\mathbb Z}$; if instead $q$ is even, 
then $k$ is odd, so $(k-1)/2\in{\mathbb Z}$. Therefore $\eta^{q^{j}-1}=1$, 
so $\eta\in\F_{q^j}$, hence $j=n$ by irreducibility of $\chi_{\x}$. 
\epf

\begin{obs}\label{punto:ss-quasi-real}
If $\St\in \SL_2(q)$ is semisimple 
irreducible such that 
$\St^q=\lambda \St$, then $q\equiv 3 \mod 4$. Indeed, by 
Remark \ref{obs:standard} \ref{punto:ss-c} its minimal polynomial is 
$X^2+1$ which is irreducible only if $q\equiv 3 \mod 4$. 
Thus, a semisimple irreducible element  $\pi(\St)$ in $\PSL_2(q)$ is 
quasi-real unless $q\equiv3 \mod 4$ and $\chi_{\St}=X^2+1$.  
\end{obs}


\subsection{Subgroups of $\PSL_2(q)$}\label{subsec:dickson}

In this subsection $\Gb = \PSL_2(q)$. We stress that we assume $q\neq 2$, $3$, $4$, $5$, $9$
to avoid coincidences with cases treated elsewhere,
see \cite[Subsection 1.2]{ACG-I}. We recall Dickson's classification of
all subgroups of $\Gb$. 
Let $d = (2, q-1)$.

\begin{theorem}\label{th:subgps-psl2}\cite[Theorems 6.25, p. 412; 6.26, p. 414]{suzuki}
A  subgroup of $\PSL_2(q)$, $q=p^m$ is isomorphic to one of the following groups.

\begin{enumerate}[leftmargin=*,label=\rm{(\alph*)}]
  \item\label{dickson:point-dihedral} The dihedral groups of order $2(q\pm 1)/d$ and their subgroups. 
  There are always such subgroups.
  \item\label{dickson:point-elementary-abelian}
  A group $H$ of order $q(q-1)/d$ and its subgroups. It has a normal $p$-Sylow subgroup $Q$ that is elementary
  abelian and the quotient $H/Q$ is cyclic of order $(q-1)/d$. There are always such subgroups.
  \item\label{dickson:point-a4} $\mathbb A_4$, and there are such subgroups except when $p=2$ and $m$ is odd.
    \item\label{dickson:point-s4} $\mathbb S_4$, and there are such subgroups if and only if $q^2 \equiv 1 \mod 16$.
      \item\label{dickson:point-a5} $\mathbb A_5$, and
      there are such subgroups if and only if $q(q^2 - 1) \equiv 0 \mod 5$.
  \item\label{dickson:point-psl2} $\PSL_{2}(t)$ for some
  $t$ such that $q = t^h$, $h\in \N$.
  There are always such subgroups.
  \item\label{dickson:point-pgl2} $\PGL_{2}(t)$ for some
  $t$ such that $q = t^h$, $h\in \N$. If $q$ is odd, then
 there are such subgroups  if and only if $h$ is even and $q = t^h$.
 Note that for $q$ even, this reduces to case \ref{dickson:point-psl2}. \qed
\end{enumerate}
\end{theorem}

For further use we record the following consequence of  Theorem \ref{th:subgps-psl2}. 
If $q$ is odd, then involutions in 
$G$ are semisimple. By looking at the possible eigenvalues of such an 
element we see that there is only one class of 
non-trivial involutions in $\Gb$ for $q$ odd. 

\begin{cor}\label{cor:sl2-involutions} Let $\Oc$ be the conjugacy class of non-trivial 
involutions in $\Gb$.
If $q > 7$ is odd, then $\Oc$ is of type D; 
while  $\Oc$ is of type C, when $q = 7$.
\end{cor}

\pf Assume that $q > 7$ is odd. Recall that the case $q=9$ is excluded. 
By Theorem \ref{th:subgps-psl2} \ref{dickson:point-dihedral}, $\Gb$ contains  
dihedral subgroups $D_1$ and $D_2$ of order $q-1$ and $q+1$ respectively. If $q\equiv 1 \mod 4$,
 resp. $q\equiv 3 \mod 4$, then  $\Oc \cap D_1$, respectively $\Oc \cap D_2$, is 
 of type D by \cite[Lemma 2.1]{AFGV-simple}.
 
Assume now that $q = 7$. Then the claim follows from Lemma \ref{lem:sl32-21-C} because 
 $\PSL_2(7)\simeq \PSL_3(2)$ and $\Oc$ is the class of involutions therein.
\epf

\subsection{Subgroups of  $\PSL_{3}(q)$}\label{subsec:psl3-ss}
In this subsection we recall some classification results about the
subgroups of $\Gb=\PSL_3(q)$ and $\PSU_{3}(q)$.
The classification of the subgroups of $\Gb$ for $q$ 
odd was obtained by Mitchell in 1911, whereas the classification
of maximal subgroups in $\Gb$ for $q$ even was achieved by Hartley in 1925.  
We set $d=(q-1,3)$.
If $q$ is even, $d=1$ 
exactly when $q$ is not a square.

\begin{theorem}\label{th:subgps-psl3-odd}\cite[Theorems 1.1, 7.1]{B}

Let $K$ be a subgroup of $\PSL_3(q)$ with $q=p^m$ and $p$ an odd prime.

\emph{(I)} Assume that $K$ has no non-trivial normal elementary abelian subgroup. 
Then $K$ is isomorphic to 
one of the following groups.

\begin{enumerate}
[leftmargin=*,label=\rm{(\alph*)}]
  \item\label{bloom:point-psl3} $\PSL_{3}(t)$ for some $t = p^a$ such that $q = t^h$, $h\in \N$.

  \item\label{bloom:point-pu3} $\PSU_{3}(t)$ for some $t = p^a$ such that $2a \mid m$.

  \item\label{bloom:point-psl3-extension} If $t = p^a$ satisfies $t \equiv 1 \mod 3$ and  $3a \mid m$,
  then there
  is a subgroup containing the subgroup of type \ref{bloom:point-psl3} with index 3.

  \item\label{bloom:point-pu3-extension} If $t = p^a$ satisfies $t \equiv 2 \mod 3$ and  $6a \mid m$, 
  then there
  is a subgroup containing the subgroup of type \ref{bloom:point-pu3} with index 3.

  \item\label{bloom:point-psl2} $\PSL_{2}(t)$ or $\PGL_{2}(t)$ for some $t = p^a \neq 3$ such that
  $q = t^h$, $h\in \N$.

  \item\label{bloom:point-psl2(5)} $\PSL_{2}(5)$, when $q \equiv \pm 1 \mod 10$.

  \item\label{bloom:point-psl2(7)} $\PSL_{2}(7)$, when $q^3 \equiv  1 \mod 7$.


   \item\label{bloom:point-a67} $\mathbb A_6$, $\mathbb A_7$ or a group containing $\mathbb A_6$
   with index 2, when $p = 5$
   and $m$ is even.

   \item\label{bloom:point-a6} $\mathbb A_6$, when $q \equiv  1$    or $19 \mod 30$.
\end{enumerate}

There are always such subgroups under the indicated restrictions.

\emph{(II)} Assume that $K$ has a non-trivial normal elementary abelian subgroup. 
Then one of the following happens:

\begin{enumerate}
[leftmargin=*,label=\rm{(\alph*)}]\setcounter{enumi}{9}
  \item\label{bloom:point-cyclic-normal} $K$ has a cyclic $p$-regular normal subgroup of index $\leq 3$.

\item\label{bloom:point-diagonal-normal} $K$ has a diagonal normal subgroup $L$ such that $K/L$ is isomorphic
to a subgroup
of $\st$.

\item\label{bloom:point-diagonal-diese} $K$ has a normal elementary-abelian p-subgroup $H$ such that $K/H$
is isomorphic
to a subgroup of $\GL_2(q)$. We include the case $H = \{1\}$.

\item\label{bloom:point-diagonal-abelian33} $K$ has a normal abelian subgroup $H$  of type $(3, 3)$, with $K/H$
isomorphic to a subgroup of $\SL_2(3)$. All subgroups of $\SL_2(3)$ do occur in this
context. This happens when $q \equiv 1 \mod 9$.

\item\label{bloom:point-diagonal-quaternion} $K$ has a normal abelian subgroup $H$  of type
$(3, 3)$, with $K/H$ isomorphic to a subgroup of the quaternion group $\Q$ of order 8.
All subgroups of $\Q$ do occur in this context. This happens when $q \equiv 1 \mod 3$, $q \not\equiv 1 \mod 9$.
\qed
\end{enumerate}

\end{theorem}

The following theorem gives the classification of the maximal subgroups of $\Gb$ for
$q$ even.

\begin{theorem}\label{th:subgps-psl3-even}\cite[Theorem 8, Summary]{hartley}
Let $M$ be a maximal subgroup  of $\PSL_3(q)$, with $q=2^m$. Then $M$ is one of the following:
\begin{enumerate}[leftmargin=*,label=\rm{(\alph*)}]
\item\label{hartley:normalizer-tori}A subgroup of order $q^3(q+1)(q-1)^2/d$ or $6(q-1)^2/d$. 
\item\label{hartley:normalizer-maximal torus}The normalizer of a  maximal torus of order
$(3)_q/d$. The torus has index $3$ in $M$.
\item\label{hartley-psl} $\PSL_3(2^k)$, where $m/k$ is prime.
\item\label{hartley-pgl} A group containing $\PSL_3(2^{2a})$ as a normal subgroup of index $3$. This happens if $m=6a$.
\item\label{hartley-psu} $\PSU_3(t)$. This happens when $q=t^2$ is square.
\item\label{hartley-pgu} A group containing $\PSU_3(2^a)$ as a normal subgroup of index $3$. This happens if 
$a$ is odd and  $m = 6a$.
\item\label{hartley-a6} A group  isomorphic to ${\mathbb A}_6$. This happens when $q = 4$.\qed
\end{enumerate}
\end{theorem}

For inductive arguments we will also need the classification of maximal subgroups for $\PSU_3(q)$, for $q$ even, 
also due to Hartley.

\begin{theorem}\label{th:subgps-psu3-even}\cite{hartley}
Let $q = 2^m$.
Let $M$ be a maximal subgroup  of $\PSU_3(q)$. Then $M$ is one of the following subgroups, where
$e:=(3, q+1)$:
\begin{enumerate}[leftmargin=*,label=\rm{(\alph*)}]
\item\label{hartley:psu-normalizer-tori}A subgroup of order $q^3(q+1)(q-1)/e$,
$q(q+1)^2(q-1)/e$ or $6(q+1)^2/e$.
\item\label{hartley:psu-normalizer-maximal torus}
The normalizer of a maximal torus of order
$3(q^2-q+1)/e$. The torus has index $3$ in $M$.
\item\label{hartley-psu-psu} $\PSU_3(2^l)$, where $m/l$ is an odd prime. 
If $l=1$ this group has order $72$.
\item\label{hartley-pgu-psu} A group containing $\PSU_3(2^l)$ as a normal subgroup of index $3$. This happens if 
$l$ is odd and  $m = 3l$. If $l=1$ this group has order $216$.
\item\label{hartley-36} A group of order $36$. This happens when $q = 4$.\qed
\end{enumerate}
\end{theorem}

\subsection{Reduction to the irreducible case}\label{subsec:reduction-irred}

In this subsection we prove that all non-irreducible semisimple conjugacy 
classes in $\Gb = \PSL_{n}(q)$  collapse. 
The proof is split in several lemmata in which we use different criteria for a rack to collapse.
Recall that $G = \SL_{n}(q)$ and $\pi: G \to \Gb$ is the natural projection.

We begin with the case where $\Oc$ is a conjugacy class of a diagonal and non-central element.

\begin{lema}\label{lema:reduction-diagonal}
If $T\in G$ is diagonal but not central, then $\Oc_{\pi(T)}^{\Gb}$ collapses.
\end{lema}

\pf Assume first that $n > 2$. Say that $T = \diag (a_1\dots, a_n)$ with $a_j \in \F_q^{\times}$; 
since $T$ is not central,
at least two of the $a_j$'s are different, so $q>2$. 
Assume that $a_1 \neq a_2$. and consider the following 
subsets of $\Oc^{G}_{T}$:
\begin{align*}
X_{1}  &= \left\{ r_c:=
\left(\begin{smallmatrix}
a_1 & c & 0 &\cdots & 0\\
0 & a_2 & 0 &\cdots & 0 \\
\vdots &\vdots & \vdots &\ddots & \vdots \\
0 & 0 & 0 &\cdots & a_n\end{smallmatrix}\right):\ c\in \F_{q}
\right\}, \\
X_{2}  &= \left\{s_f:=
\left(\begin{smallmatrix}
a_2 & f & 0 &\cdots & 0\\
0 & a_1 & 0 &\cdots & 0 \\
\vdots &\vdots & \vdots &\ddots & \vdots \\
0 & 0 & 0 &\cdots & a_n\end{smallmatrix}\right):\ f\in \F_{q}
\right\}.
\end{align*}

Since \begin{align}\label{eq:reduction-case1}
\left(\begin{smallmatrix}
d & c\\
0 & e
\end{smallmatrix}\right) \trid \left(\begin{smallmatrix}
e & f\\
0 & d
\end{smallmatrix}\right) =
\left(\begin{smallmatrix}
e & de^{-1}(c + f)-c\\
0 & d
\end{smallmatrix}\right),
\end{align}
$Y = X_{1} \coprod X_{2}$ is a decomposable subrack  
of $\Oc^{\SL_{n}(q)}_{T}$, cf. \eqref{eq:orbit-ss}.
Set $ r = r_1$ and $s = s_0$, so that $r \trid s \neq s$. 
Then $r_c \trid s = s_{(a_1a_2^{-1}-1)c}$, hence $Y\trid s = X_{2}
=\Oc_s^{\langle Y\rangle}$ by \eqref{eq:reduction-case1};
similarly, $Y\trid r = X_{1}=\Oc_r^{\langle Y\rangle}$. 
Also $\vert X_1\vert =\vert X_2\vert = q > 2$.
Hence $\Oc^{G}_{T}$ is of type C and 
since $\pi_{\vert Y}$ is injective,  $\Oc^{\Gb}_{\pi(T)}$ is of type C. 

\smallbreak
Assume next that $n = 2$, so that $T=\left(\begin{smallmatrix}
      a & 0\\
      0 & a^{-1}
     \end{smallmatrix}\right)$, with $a^{2}\neq 1$. Let
\begin{align*}
X_{a}  &= \left\{r_c :=
\left(\begin{smallmatrix}
a & c\\
0 & a^{-1}
\end{smallmatrix}\right):\ c\in \F_{q}
\right\}, &
X_{a^{-1}}  &= \left\{s_f :=
\left(\begin{smallmatrix}
a^{-1} & f\\
0 & a
\end{smallmatrix}\right):\ f\in \F_{q}
\right\}.
\end{align*}
Then $Y = X_{a} \coprod X_{a^{-1}}$ is a
decomposable subrack  of $\Oc^{\SL_{2}(q)}_{T}$,  by \eqref{eq:reduction-case1}; and 
$r_c \trid s_f \neq s_f$ if and only if $c+f \neq 0$; respectively,
$r_c\trid(s_f\trid(r_c\trid s_f)) \neq s_f$ if and only if $2(c+f) \neq 0$.

\smallbreak
\emph{Assume that $q$ is odd}. Then $\Oc^{\SL_{2}(q)}_{T}$ is of type D. 
If $a^{4}\neq 1$, then $\pi_{\vert Y}$ is injective, 
hence $\Oc^{\Gb}_{\pi(T)}$ is of type D. 
If $a^{4} = 1$, then $\pi(T)$ is an involution and $q\neq 7$, 
hence $\Oc^{\Gb}_{\pi(T)}$ is of type D
by Corollary \ref{cor:sl2-involutions}.

\smallbreak
\emph{Assume that  $q$ is even}. Set $ r = r_1$ and $s = s_0$; thus $r \trid s \neq s$. 
Now $r_c \trid s = s_{c+ca^2}$ by \eqref{eq:reduction-case1}, 
hence $Y\trid s = X_{a^{-1}}$;
similarly, $Y\trid r = X_{a}$. Also $\vert X_{a^{\pm 1}}\vert = q > 2$.
Hence $\Oc^{G}_{T}$ is of type C, and the claim follows since  $\Gb = G$.
\epf

\begin{lema}\label{lema:reduction-block-eigenvalue}
Let $T \in \SL_{n}(q)$ semisimple not diagonal,
with at least one eigenvalue.
Then $\Oc_{\pi(T)}^{\Gb}$ is of type C.
\end{lema}
\pf
By hypothesis, we may assume that either $T=\left(\begin{smallmatrix}
 a & 0 & 0\\
 0 & B &0 \\ 0 & 0 & C
\end{smallmatrix}\right)$, where $B \in \GL_{e}(q)$ is semisimple irreducible not diagonal, 
$C$ is semisimple and $a = (\det B \det C)^{-1}$, or else $T=\left(\begin{smallmatrix}
 a & 0 \\
 0 & B  
\end{smallmatrix}\right)$, where $B \in \GL_{e}(q)$,  
is semisimple irreducible not diagonal, $e = n-1$,
and $a = (\det B)^{-1}$. 
In both cases $B\neq B^{q}\in \Oc_{B}^{\GL_{e}(q)}$ 
and $\det B^{q} = (\det B)^{q} = \det B$.
We treat only the first possibility, the second being analogous.
Let
\begin{align*}
X_{1} &= \left\{ r_v:=
\left(\begin{smallmatrix}
a & v & 0\\
 0 & B &0 \\ 0 & 0 & C
\end{smallmatrix}\right):v\in \F_{q}^{e}
\right\} , &
X_{2}  &= \left\{ s_w :=
\left(\begin{smallmatrix}
a & w & 0\\
 0 & B^q &0 \\ 0 & 0 & C 
\end{smallmatrix}\right):  w \in \F_{q}^{e}
\right\}.
\end{align*}
Then $r_v\trid s_w = s_{\left(a(w-v) + v B^q \right)B^{-1}}$; hence 
$Y = X_{1} \coprod X_{2}$ is a decomposable subrack  
of $\Oc^{\SL_{n}(q)}_{T}$, cf. \eqref{eq:orbit-ss}. 
If $s= s_0$, then $r_v\trid s = s_{v\left(-a +  B^q \right)B^{-1}}$, thus $Y\trid s = X_2$
because $-a +  B^q$ is invertible by hypothesis. 
Similarly, $Y\trid r = X_1$. Also, if $v \neq 0$ and $r= r_v$, then  
$r\trid s \neq s$. Since $\vert X_1\vert =  \vert X_2\vert = q^e > 2$, $\Oc_{T}^{G}$
is of type C. Since $\pi_{\vert Y}$ is injective, $\Oc_{\pi(T)}^{\Gb}$ is of type C.
\epf

Next we treat the cases where the non-irreducible semisimple element has no 
eigenalues in $\F_{q}$.

\begin{lema}\label{lema:reduction>2blocks}
Let $T \in \SL_{n}(q)$ semisimple not irreducible, 
with at least $3$ irreducible blocks and no eigenvalues in $\F_q$.
Then $\Oc_{\pi(T)}^{\Gb}$ is of type C.
\end{lema}
\pf We may assume that $T=\left(\begin{smallmatrix}
 A & 0 & 0\\
 0 & B & 0\\
 0 & 0 & C
\end{smallmatrix}\right)$, where $A \in \GL_{d}(q)$, $B \in \GL_{e}(q)$ and 
$C \in \GL_{f}(q)$ are semisimple not diagonal, $A$ and $B$ are irreducible, $d+e+f=n$
and $\det A \det B \det C=1$. 

As in the previous proof,
$B\neq B^{q}\in \Oc_{B}^{\SL_{e}(q)}$ and $\det B^{q} = \det B$.
Let 
\begin{align*}
R &= \left\{\left(\begin{smallmatrix} E & 0 &0 \\ 
 0& B & 0\\
 0 & 0 &C\end{smallmatrix}\right):
\ E\in \Oc_{A}^{\SL_{d}(q)}\right\}, &
S &= \left\{\left(\begin{smallmatrix} F & 0 &0 \\ 
 0& B^q & 0\\
 0 & 0 &C\end{smallmatrix}\right):
\ F\in \Oc_{A}^{\SL_{d}(q)}\right\},
\end{align*}
Then
$Y=R\coprod S$ is a decomposable subrack of $\Oc_{T}^{\SL_{n}(q)}$. 
A normal subgroup of $\GL_{d}(q)$ is either central or contains $\SL_{d}(q)$
\cite[p. 40]{D}. Let $K = \langle \Oc_{A}^{\SL_{d}(q)}\rangle \nsbgp \GL_{d}(q)$.
Since $A$ is not central, $\SL_{d}(q) \le K$ and $K$ is not abelian.
Therefore there are
$x,y$ in $\Oc_{A}^{\GL_{d}(q)}$
such that $x\trid y \neq y$. By \eqref{eq:orbit-ss}, 
$\Oc_{A}^{\GL_{d}(q)}=\Oc_{A}^{\SL_{d}(q)}$. Then  $r:=\left(\begin{smallmatrix}
 x & 0 & 0\\
 0 & B & 0 \\
 0 & 0 & C
\end{smallmatrix}\right)$ and $s :=\left(\begin{smallmatrix}
 y & 0 &0\\
 0 & B & 0\\
 0 & 0 & C
\end{smallmatrix}\right)$ satisfy $r\trid s \neq s$. 
By the same reason, $\vert R\vert = \vert S\vert = \vert \Oc_{T}^{\SL_{d}(q)}\vert > 2$.
Finally, $ \GL_{d}(q) \trid y \supseteq K \trid y \supseteq \SL_{d}(q) \trid y = \GL_{d}(q) \trid y$.
Hence \eqref{eq:typeC-rack-indecomposable} holds and $\Oc_{T}^{G}$
is of type C. Since $\pi_{\vert Y}$ is injective, $\Oc_{\pi(T)}^{\Gb}$ is of type C.
\epf

Now we prove that the 
conjugacy class of a non-irreducible non-diagonal semisimple element with 
two irreducible blocks and no eigenvalues in $\F_{q}$ collapses. Again, we split 
the proof in several lemmata, depending on the relation between 
the blocks, $n$ and $q$.

\begin{lema}\label{lema:reduction-2blocks}
Let $T \in \SL_{n}(q)=\left(\begin{smallmatrix}
 A & 0 \\
 0 & B  
\end{smallmatrix}\right)$ semisimple, with $A\in GL_d(q)$ and $B\in\GL_e(q)$ 
irreducible and without eigenvalues in $\F_q$. 
If $B^q\not \in Z(\SL_e(q))B$, or else  $A^q\not \in Z(\SL_d(q))A$,
then $\Oc_{\pi(T)}^{\Gb}$ is of type C.
\end{lema}
\pf Up to interchanging the role of $A$ and $B$ we may 
assume that $B^q\not \in Z(\SL_e(q))B$.
Arguing as in the previous proof, we take 
\begin{align*}
R &= \left\{\left(\begin{smallmatrix} C & 0 \\ 0& B\end{smallmatrix}\right):
\ C\in \Oc_{A}^{\SL_{d}(q)}\right\}, &
S &=\left\{\left(\begin{smallmatrix} D & 0 \\ 0& B^{q}\end{smallmatrix}\right):
\ D\in \Oc_{A}^{\SL_{d}(q)} \right\}.
\end{align*}

Then $Y=R\coprod S$ is a decomposable subrack of $\Oc_{T}^{\SL_{n}(q)}$ and 
for $x,y$ in $\Oc_{A}^{\SL_{d}(q)}$ such that $x\trid y \neq y$, the elements 
$r:=\left(\begin{smallmatrix}
 x & 0 \\
 0 & B  
\end{smallmatrix}\right)$ and $s :=\left(\begin{smallmatrix}
 y & 0 \\
 0 & B^q  
\end{smallmatrix}\right)$ satisfy $r\trid s \neq s$. 
The same argument as in the proof of Lemma \ref{lema:reduction>2blocks} shows that
$\vert R\vert = \vert S\vert = \vert \Oc_{A}^{\SL_{d}(q)}\vert > 2$ and that
\eqref{eq:typeC-rack-indecomposable} holds for $\Oc_{T}^{G}$. We prove that
$\pi_{\vert Y}$ is injective. Clearly $\pi_{\vert R}$ and $\pi_{\vert S}$ are injective. 
If 
$\pi\left(\begin{smallmatrix} C & 0 \\ 0& B\end{smallmatrix}\right)=
\pi\left(\begin{smallmatrix} D & 0 \\ 0& B^q\end{smallmatrix}\right)$ 
for some $C, D\in \Oc_A^{\SL_d(q)}$,
then there would be a $\lambda\in\F_q$ such that $B^q=\lambda B$. 
Computing the determinant we get
$\lambda^{e}=1$ and therefore $\lambda\id\in Z(\SL_e(q))$ contradicting our hypothesis. 
\epf

\begin{lema}\label{lema:reduction-2equal-blocks}
Let $T, A, B, d, e$ be as in Lemma \ref{lema:reduction-2blocks}
with $A$ and $B$ irreducible and without eigenvalues in $\F_q$ and let $l:=(q-1,d,e)$.
If $l\neq d$  or $l\neq e$, then $\Oc_{\pi(T)}^{\Gb}$ is of type C.
\end{lema}
\pf Up to interchanging the role of $A$ and $B$ we may assume that $l\neq e$.
Arguing as  in the proof of Lemma \ref{lema:reduction-2blocks}, 
with same $Y$, we conclude that 
$\Oc_{T}^{\SL_{n}(q)}$ is of type C.
We show that $\pi_{\vert Y}$ is injective.
If  $\pi\left(\begin{smallmatrix} C & 0 \\ 0& B\end{smallmatrix}\right)=
\pi\left(\begin{smallmatrix} D & 0 \\ 0& B^q\end{smallmatrix}\right)$ 
for some $C, D\in \Oc_A^{\SL_d(q)}$, then
there would be $\lambda\in\F_q$ such that $B^q=\lambda B$, 
$D=\lambda C$. Computing the determinants we get
$\lambda^{d}=\lambda^e=1$ and therefore $\lambda^{l}=1$. 
By Remark \ref{obs:standard} \ref{obs:Cq-scalar}, $\lambda$
would be  a primitive $e$-th root of $1$, contradicting our hypothesis.
\epf

\begin{lema}\label{lema:reduction-2equal-blocks-not-conjugate}
Let $T \in \SL_{n}(q)=\left(\begin{smallmatrix}
 A & 0 \\
 0 & B  
\end{smallmatrix}\right)$ be semisimple, with $A, B\in \GL_d(q)$ 
irreducible with no eigenvalues in $\F_q$.
If $B\not\in Z(\SL_d(q))\Oc_A^{\SL_d(q)}$, then 
$\Oc_{\pi(T)}^{\Gb}$ is of type C.
\end{lema}
\pf Since $A$ and $B$ are irreducible, then each of them  lies in 
a maximal torus of $\GL_d(q)$
of order $q^d-1$. 
All such tori are conjugate in $\GL_d(q)$ \cite[Proposition 25.1]{MT}. 
So, up to replacing $B$ by   $B'\in \Oc_B^{\GL_d(q)}=
\Oc_B^{\SL_d(q)}$ in the same torus as $A$,
we can assume that $A$ and $B$ commute. 
We set
\begin{align*}
R &= \left\{\left(\begin{smallmatrix} D & 0 \\ 0& B\end{smallmatrix}\right):
\ D\in \Oc_{A}^{\SL_{d}(q)}\right\}, &
S &=\left\{\left(\begin{smallmatrix} E & 0 \\ 0& A\end{smallmatrix}\right):
\ E\in \Oc_{B}^{\SL_{d}(q)} \right\}.
\end{align*}
Then $Y=R\coprod S$ is a decomposable subrack of $\Oc_{T}^{\GL_{n}(q)}$.
The same argument as above gives $\vert R\vert , \vert S\vert  > 2$
and \eqref{eq:typeC-rack-indecomposable} for $\Oc_{T}^{G}$. 
Therefore $\Oc_{T}^{G}$ is of type C if we can find 
noncommuting  $D\in \Oc_{A}^{\SL_{d}(q)}$ and $E\in \Oc_{B}^{\SL_{d}(q)}$. 
Let $D\in \Oc_{A}^{\SL_{d}(q)}$ be in a  maximal torus $S_D$. 
Since $D$ is irreducible, it has distinct eigenvalues so  
its centralizer is $S_D$. Therefore it is enough to choose $E\not\in S_D$. 
If $B\not\in S_D$ we take $E=B$. Assume $B\in S_D$ and let
 $x\in \SL_d(q)$ such that $x\not\in N_{\SL_d(q)}(S_D)$. Then 
$E:=xB x^{-1}\in x S_D x^{-1}\neq S_D$. As every semisimple element 
with distinct eigenvalues lies in a unique maximal torus \cite[2.3]{Hu}, $E\not\in S_D$. 
Injectivity of $\pi_{\vert Y}$ yields the statement. 
\epf

We are left with the analysis of $\Oc_{\pi(T)}^{\Gb}$ for 
$T \in \SL_{n}(q)=\left(\begin{smallmatrix}
 A & 0 \\
 0 & B  
\end{smallmatrix}\right)$ semisimple, with $A, B\in \GL_d(q)$ 
irreducible, $B\in Z(\SL_d(q)) \Oc_A^{\GL_d(q)}$,  
$A^q\in Z(\SL_d(q))A$, $\det(A)\det(B)=1$, $d>1$ and $d\mid(q-1)$. 
It follows from Remark \ref{obs:standard} \ref{obs:Cq-scalar} 
that under these assumptions 
$Z(\SL_d(q)) \Oc_A^{\GL_d(q)} =\Oc_A^{\GL_d(q)}$ so 
$B\in \Oc_A^{\GL_d(q)}$, $\det(A)=\det(B)=\pm1$ 
and the characteristic polynomial of $A$ is $X^d+(-1)^d\det(A)$. 
This polynomial is irreducible only if $d$ is even, $\det(A)=1$ 
and $-1$ is not a square in $\F_q$, i.e., 
$q\equiv 3 \mod 4$.
We analyze this situation, studying separately the cases $d>2$  and $d=2$. 

\begin{lema}\label{lema:reduction-2equal-blocks-conjugate}
Let $T \in \SL_{2d}(q)=\left(\begin{smallmatrix}
 A & 0 \\
 0 & A  
\end{smallmatrix}\right)$ be semisimple, with $A\in \SL_d(q)$ 
irreducible, $d>2$,   $A^q=\mu A$,  
$d|q-1$, $d$ even and $\mu$ a primitive $d$-th root of $1$. Then 
$\Oc_{\pi(T)}^{\Gb}$ is of type D.
\end{lema}
\pf It is always possible to find $a\neq b\in \I_{d-1}$ such that, 
$\mu^{2(a+b)}\neq 1$.  The matrices 
$x=\diag(A, \mu^a A)$ and $y=\diag (A, \mu^b A)$ lie in $\Oc_T^{\SL_{2d}(q)}$. 
We set
\begin{align*}
R &= \left\{X=\left(\begin{smallmatrix} A & X' \\ 
 0& \mu^a A\end{smallmatrix}\right):
\ X\in \Oc_{T}^{\SL_{2d}(q)}\right\}, \\
S &= \left\{Z=\left(\begin{smallmatrix}
 A & Z' \\ 
 0& \mu^b A
\end{smallmatrix}\right):
\ Z\in \Oc_{T}^{\SL_{2d}(q)}\right\}.
\end{align*}
Then $Y:=R\coprod S$ is a decomposable subrack. 
Let 
\begin{align*}r=\left(\begin{smallmatrix} \id_d & \id_d \\ 
 0& \id_d\end{smallmatrix}\right)\trid x=\left(\begin{smallmatrix} A & (\mu^a-1)A \\ 
 0& \mu^a A\end{smallmatrix}\right)\in R\end{align*} and $s:=y\in S$. 
A direct computation shows that $(rs)^2\neq (sr)^2$. Since for our choice of $a$ and $b$ the map $\pi_{\vert Y}$ is injective, 
$\Oc_{\pi(T)}^{\Gb}$ is of type D.
\epf

\begin{lema}\label{lema:reduction-2equal-blocks-2}
Let $T \in \SL_{4}(q)=\left(\begin{smallmatrix}
 A & 0 \\
 0 & A  
\end{smallmatrix}\right)$ be semisimple, with $A\in \GL_2(q)$ irreducible, $A^q=- A$,  
$q\neq 3$. Then $\Oc_{\pi(T)}^{\Gb}$ is of type D.
\end{lema}
\pf By Remark \ref{obs:standard} \ref{obs:Cq-scalar}, 
the characteristic polynomial of $A$ is necessarily $X^2+1$ so
$q\equiv 3 \mod 4$, and $\pi(T)$ is an involution.  
Let $\zeta$ be a generator of $\F_q^\times$.
By \cite[Lemma 2.5]{FaV}, $\Oc_{\pi(T)}^{\Gb}$ 
is of type D if and only if there exist $r,s\in \Oc_{\pi(T)}^{\Gb}$
such that $\ord (rs)>4$ is even. 
Let 
\begin{align*}\erre:=\left(\begin{smallmatrix}
0&-\zeta&0&0\\
\zeta^{-1}&0&0&0\\
0&0&0&1\\
0&0&-1&0
\end{smallmatrix}\right),&&
\esse:=\left(\begin{smallmatrix}
0&0&0&1\\
0&0&1&0\\
0&-1&0&0\\
-1&0&0&0
\end{smallmatrix}\right).\end{align*}
As $\erre$ and $\esse$ are semisimple matrices 
with characteristic polynomial equal to $(X^2+1)^2$, 
they lie in $\Oc_T^{\SL_4(q)}$. In addition,
\begin{align*}\erre\esse=\left(\begin{smallmatrix}
0&0&-\zeta&0\\
0&0&0&\zeta^{-1}\\
-1&0&0&0\\
0&1&0&0
\end{smallmatrix}\right),&&(\erre\esse)^2=\diag(\zeta,\zeta^{-1},\zeta,\zeta^{-1}),\end{align*}
so, for $r:=\pi(\erre)$, $s:=\pi(\esse)$ we have 
$\ord(rs)=2\ord(\zeta^2)=q-1>4$, as $q\neq 3$, 
$q\equiv 3 \mod4$. 
\epf

\begin{lema}\label{lema:reduction-d2-q3}
Let $T \in \SL_{4}(3)=\left(\begin{smallmatrix}
 A & 0 \\
 0 & A  
\end{smallmatrix}\right)$ be semisimple, with $A=\left(\begin{smallmatrix}
 0& 2 \\
 1 & 0  
\end{smallmatrix}\right)$. Then $\Oc_{\pi(T)}^{\Gb}$ is of type D.
\end{lema}
\pf  Let 
\begin{align*}v:=\left(\begin{smallmatrix}
1&0&0&1\\
1&1&2&1\\
1&1&0&0\\
0&0&0&1
\end{smallmatrix}\right),&&
\esse:=v\trid T=\left(\begin{smallmatrix}
1&2&0&0\\
2&2&0&0\\
2&0&2&1\\
0&2&1&1
\end{smallmatrix}\right)\in \Oc_{\pi(T)}^{\SL_4(q)}.\end{align*}
A direct computation shows that $|T \esse| =12$, so setting 
$r:=\pi(T)$ and $s:=\pi(\esse)$ we have 
$\ord(rs)\in \{6,12\}$. We apply \cite[Lemma 2.5]{FaV}.
\epf

\medskip

We end this subsection with the statement of the following result.

\begin{prop}\label{prop:reduction}
If $T\in \SL_{n}(q)$ is semisimple and not irreducible, 
then $\Oc_{T}^{\PSL_{n}(q)}$ collapses. \qed
\end{prop}

\pf Let $T\in \SL_{n}(q)$ be a non-irreducible semisimple element in $\PSL_{n}(q)$. 
If $T$ is has at least one eigenvalue and it is not central, the claim follows by 
Lemmata \ref{lema:reduction-diagonal}, \ref{lema:reduction-block-eigenvalue}. Assume
$T$ has no eigenvalues in $\F_{q}$. If $T$ has
at least 3 irreducible blocks, we apply Lemma 
\ref{lema:reduction>2blocks}, and if $T$ has exactly 
2 irreducible blocks, the assertion follows by 
Lemmata \ref{lema:reduction-2blocks}, \ref{lema:reduction-2equal-blocks},
\ref{lema:reduction-2equal-blocks-not-conjugate} 
\ref{lema:reduction-2equal-blocks-conjugate}, 
\ref{lema:reduction-2equal-blocks-2} and \ref{lema:reduction-d2-q3}. 
\epf

\subsection{General results on irreducible semisimple classes}
\label{subsec:general-irreducible}
Set $\Gb=\PSL_n(q)$. In this subsection we prove that for certain subgroups 
$K\leq \Gb$ and $x\in \Gb$ semisimple, 
we have that $\Oc_{x}^{\Gb} \cap K  = \Oc^{K}_{x}$.
These results will be used in the sequel. 

\begin{lema}\label{lem:one-gcd-in-pslt}
Let $K=\PSL_n(t)\leq  \Gb$ for $t=p^a$ and $a|m$ with $(t-1,n)=(q-1,n)$.
If $x\in K$ is semisimple, then $\Oc_x^{\Gb}\cap K=\Oc_x^K$.
\end{lema}
\pf Let $\x\in \SL_n(t)$ such that $\pi(\x)=x$. By Remark \ref{obs:standard} \ref{punto:ss-d},
$\Oc_{\x}^{\SL_n(\kk)}\cap \SL_n(t)=\Oc_{\x}^{\SL_n(t)}$. 
A fortiori, $\Oc_{\x}^{\SL_n(q)}\cap \SL_n(t)=\Oc_{\x}^{\SL_n(t)}$.
Let $y=\pi(\y)\in \Oc_x^{\Gb}\cap K$, with $\y\in\SL_n(t)$. 
Since $(t-1,n)=(q-1,n)$, we have that $Z(\SL_n(q))=Z(\SL_n(t))$ and
then, for some $z\in Z(\SL_n(q))$ there holds 
$\y\in z(\Oc_{\x}^{\SL_n(q)}\cap \SL_n(t))=z\Oc_{\x}^{\SL_n(t)}$, so $y\in \Oc_x^K$. 
\epf

\begin{lema}\label{lem:irred-one-in-pslt}
Let $K=\PSL_n(t)\leq\Gb$ for $t=p^a$ and $a|m$. Assume in addition that 
$\left((n)_t,n\right )=1$. If $x\in K$ is semisimple irreducible, then 
$\Oc_x^{\Gb}\cap K=\Oc_x^K$.
\end{lema}
\pf Let $\x\in \SL_n(t)$ such that $\pi(\x)=x$ and let $\y\in \SL_n(t)$ with 
$y=\pi(\y)\in \Oc_{\x}^{\Gb}\cap K$. So $\y\in (z\Oc_{\x}^{\SL_n(q)})\cap \SL_n(t)$
 for some $z\id\in Z(\SL_n(q))$. Therefore it is enough to prove that $z=1$.  
Since $\ord \x$ and $\ord \y=\ord z\x$ divide $(n)_t$, we have $z^{(n)_t}=1$, 
whence the statement.
\epf

\begin{lema}\label{lem:one-in-psu} Let $K=\PSU_n(t)\leq \Gb$
for $t=p^a$ and $2a|m$. 
If $(t+1,n)=(q-1,n)$, then $\Oc_x^{\Gb}\cap K=\Oc_x^K$ 
for every semisimple $x\in K$. 
\end{lema}
\pf Let $y=\pi(\y)\in \Oc_x^{\PSL_n(q)}\cap K$, 
with $\y\in \SU_n(t)$. For some $z\in Z(\SL_n(q))$ and $\x\in \SU_n(t)$ 
with $\pi(\x)=x$ we have 
$\y\in \Oc_{z\x}^{\SL_n(q)}\cap \SU_n(t)$. Under our assumptions 
$Z(\SL_n(q))\simeq \Z/ (q-1,n)=\Z/ (t+1,n)\simeq Z(\SU_n(t))$, 
so $z\x\in \SU_n(t)$ and it is semisimple. 
The centralizer of $\x$ in $\SL_n(\kk)$ is connected, so
$\Oc_{z\x}^{\SL_n(\kk)}\cap \SU_n(t)=\Oc_{z\x}^{\SU_n(t)}$.
Thus, $\y\in  \Oc_{z\x}^{\SU_n(t)}$, and $y\in  \Oc_{x}^K$.
\epf


\section{Finite-dimensional pointed Hopf algebras over $\PSL_2(q)$}\label{sec:psl2}

\subsection{Abelian techniques}
Let $G$ be a finite group, $\oc$ a conjugacy class of $G$, $g\in
\oc$  and $(\rho,V)\in \Irr C_{G}(g)$. The abelian techniques 
are those used to conclude that
 $\dim\mathfrak{B}(\oc,\rho) = \infty$ from the consideration 
of abelian subracks and via the classification
of braided  vector spaces of diagonal type with finite-dimensional 
Nichols algebra \cite{H-all}.

\begin{lema}
\label{lem:ab-techniques-real} Assume that $\dim\mathfrak{B}(\oc,\rho)<\infty$.

\begin{enumerate}[leftmargin=*,label=\rm{(\alph*)}]
\item \cite{AZ}  If $g$ is real,  then
$\rho(g)=-1$. In particular,  $\ord g$ is even.

\item \cite{AF2,FGV1} If $g$ is quasi-real,
with $j\in \Z$ such that $g\ne g^{j}\in\oc$, then:

\renewcommand{\theenumi}{\roman{enumi}}\renewcommand{\labelenumi}{(\theenumi)}
\begin{enumerate}[leftmargin=*,label=\rm{(\roman*)}]
\item If $\deg\rho>1$, then $\rho(g)=-1$ and $g$ has even order.

\item If $\deg\rho=1$, then $\rho(g)=-1$ and $g$ has even 
order or $\rho(g)\in\mathbb{G}'_{3}$.

\item If $g^{j^{2}}\ne g$, then $\rho(g)=-1$. \qed
\end{enumerate}
\end{enumerate}
\end{lema}

\subsection{Semisimple classes in $\PSL_{2}(q)$}\label{subsec:psl2-ss}

 We recall some basic facts.

\medbreak
$\bullet$ There exist two (conjugacy classes of) maximal tori in $\SL_{2}(q)$: 
the split torus $T_1 = \{\left(\begin{smallmatrix}
      a & 0\\
      0 & a^{-1}
     \end{smallmatrix}\right), a\in \F_q^\times\}$ of order $q-1$ 
and the non-split torus $T_2$ of order 
     $q + 1$. Every non-central $x \in \SL_{2}(q)$ semisimple is
     conjugated to an element of either $T_1$ or $T_2$. Two elements $x, y\in T_i$ 
     are conjugated if and only if $x = y ^{\pm 1}$. Both $T_1$ and $T_2$ are cyclic.
		So this is the situation for $\PSL_{2}(q) = \SL_{2}(q)$ when $q$ is even.
		For uniformity of the notation, we set 
		$\Tb_i := T_i$, $i\in \I_2$, when $q$ is even.

\medbreak
 $\bullet$  Suppose $q$ is odd.
  There exist two conjugacy classes of maximal tori in $\Gb$, namely the images 
  $\Tb_1$ of the split torus $T_1$, of order $\frac{q-1}2$;
  and $\Tb_2$ of the non-split torus $T_2$, of order $\frac{q+1}2$. 
  Every $x \in \Gb$ semisimple is conjugated to an element of either $\Tb_1$ or $\Tb_2$.
Two elements $x, y\in \Tb_i$ are conjugated if and only if $x = y ^{\pm 1}$.
As a consequence, there is exactly one conjugacy class of
involutions, at most one semisimple conjugacy class of elements of order 3,
at most  one conjugacy class of elements of order 4, etc.

\begin{prop}\label{prop:sl2-ss}
Let $\Oc$ be a semisimple conjugacy class in $\PSL_2(q)$.
If $\Oc$ is not listed in Table \ref{tab:ss-psl2}, then it collapses.
\end{prop}

\pf If either $x$ is conjugated to an element in $\Tb_1$ 
or else $x$ is an involution, 
then  $\oc_x^{\Gb}$ collapses by
Lemma \ref{lema:reduction-diagonal} and Corollary \ref{cor:sl2-involutions} 
(there are no semisimple involutions when $q$ is even). 
Assume in the rest of the proof that $x\in \Gb$ is conjugated to an
element in $\Tb_2$ and $\ord x >2$.

Suppose first that $\ord x = 3$. Then, necessarily $q\equiv 2 \mod 3$. 
If moreover $q$ is odd or a square, then $\Gb$ 
contains a subgroup isomorphic to $\mathbb{A}_{4}$, by 
Theorem \ref{th:subgps-psl2}, case
\ref{dickson:point-a4}. 
Since $\Gb$
contains only one conjugacy class of elements of order 3, we 
have $Y:=\Oc_{x}^{\Gb} \cap \mathbb{A}_{4}
= \Oc_{(123)}^{\mathbb{A}_{4}}\coprod
\Oc_{(132)}^{\mathbb{A}_{4}}=R\coprod S$.  $Y$ is the so-called  \textit{cube rack},
and $R$ ans $S$ are so-called \textit{tetrahedral racks}.  They have
size 4 and do not commute with each other. Since the tetrahedral rack is indecomposable,
 $\Oc_{x}^{\Gb}$ is of type $C$ and consequently it collapses.

Now let $X$ be a subrack of $\oc^{\Gb}_x$ and $K =\langle
X\rangle$;  $X$ is a union of $K$-orbits.  We divide the 
proof with respect to the classification given in Theorem \ref{th:subgps-psl2}.

If $K$ is as in case
\ref{dickson:point-dihedral},
then $K \le \mathbb D_{\frac{q+1}{d}}$ and $X$ is abelian. 
Clearly, $K$ could not be as in case
\ref{dickson:point-elementary-abelian} because $\ord x\nmid$ the
order of such group.

Suppose that
 $K = \PSL_{2}(t)$ for some $t$ such
that $q = t^h$, $h\in \N$, case \ref{dickson:point-psl2}. Assume first that $t \neq 2,3$.  Let
$y, z\in X$. Then $\oc^K_y$ does not intersect any split torus of
$\PSL_{2}(t)$, since otherwise it would intersect $\Tb_{1}$. Let
$\mathfrak T$ be a non-split torus of $\PSL_{2}(t)$; then  $\oc^K_y$ intersects $\mathfrak T$.
Also we may
assume that either $\mathfrak T \subset \Tb_2$ or $\mathfrak T
\subset \Tb_1$ (only when $h$ is even). In the first possibility,
\begin{align*}
\emptyset \neq \oc^K_y \cap\mathfrak T \overset{\lozenge}= 
\oc^{\Gb}_y \cap \Tb_2 = \oc^{\Gb}_z \cap \Tb_2
\overset{\lozenge}= \oc^K_z \cap\mathfrak T.
\end{align*}
Here in $\lozenge$ we use that a conjugacy class intersects a torus in 
$\{x^{\pm 1}\}$. Hence $\oc^K_y= \oc^K_z = X$ is
indecomposable by \cite[Lemmata 1.9 \& 1.15]{AG-adv}. 
The second possibility is analogous. 
Now assume that $K = \PSL_{2}(2)={\mathbb S}_3$. 
Then $\Oc \cap K$ consists of $3$-cycles, hence it is abelian.
Finally, $K = \PSL_{2}(3)={\mathbb A}_4$ is excluded 
because $\Oc \cap K$ consists of involutions.

Assume that $K = \PGL_{2}(t)$ for some $t$ such that $q = t^h$, 
case \ref{dickson:point-pgl2}; we may suppose that $q$ odd,
and then $h = 2k \in \N$ should be even. If $x\in K$ is semisimple, 
then $\ord x$ divides $\vert \PGL_{2}(t)\vert = t (t^2 -1)$,
hence $\ord x$ does not divide $q + 1 = t^{2k} + 1$ and $x$ is split 
in $\PSL_{2}(q)$. In other words, $X\cap K = \emptyset$.

Hence, the only possible cases where
$X$ might not be sober are when $K\simeq \mathbb{A}_{4}$, $\mathbb{S}_{4}$
or $\mathbb{S}_{5}$, that is cases \ref{dickson:point-a4}, \ref{dickson:point-s4}
and \ref{dickson:point-a5}, respectively.
From the previous considerations, the next statement follows at once.

\begin{stepsldos}\label{step:sl2-seis}
If $\ord x > 5$, then $\oc_x^{\Gb}$ is sober.
\end{stepsldos}

We next analyze the low order cases.

\begin{stepsldos}\label{step:sl2-tres}
If $\ord x = 3$ and $q$ is even and not a square, $X$ is sober. 
\end{stepsldos}

If $q$ is even and not a square,
then cases \ref{dickson:point-a4} and
\ref{dickson:point-s4} are not possible since
$\Gb$ does not contain a subgroup isomorphic to 
$\mathbb{A}_{4}$. If $K \simeq \aco$, case \ref{dickson:point-a5}, 
then $X = \oc^{\aco}_{3}$ is
indecomposable.
 
\begin{stepsldos}\label{step:sl2-cuatro}
If $\ord x = 4$, then $\oc_x^{\Gb}$ is sober.
\end{stepsldos}

Cases \ref{dickson:point-a4} and \ref{dickson:point-a5}
are clearly excluded. If $K \simeq \sk$, case
\ref{dickson:point-s4}, then $X = \oc^4_4$ is indecomposable by 
\cite[Lemmata 1.9,1.15]{AG-adv} 
(as $\oc^4_4$ generates $\sk$ which is centerless).

\begin{stepsldos}\label{step:sl2-cinco}
The conjugacy classes of elements of order 5 in ${\Gb}$ are sober.
\end{stepsldos}

Cases \ref{dickson:point-a4} and \ref{dickson:point-s4}
are clearly excluded. Assume that $K \simeq \aco$, case
\ref{dickson:point-a5}. There are two conjugacy classes $\oc_1$
and $\oc_2$ of elements of order 5 in $\aco$ and each of them is
real, i. e. stable under inversion. Suppose that $\oc_1\subset X$
and pick $x\in \oc_1$. Then $\oc_x^{\Gb} = \oc_{x^4}^{\Gb}$ and
$\oc_{x^2}^{\Gb} = \oc_{x^3}^{\Gb}$ are the two different conjugacy
classes  of elements of order 5 in $\Gb$. If $X =
\oc_1\coprod\oc_2$, then $x^2$ would belong to $\oc_x^{\Gb}$, a
contradiction. Hence $X = \oc_1$ is indecomposable.
\epf

\subsection{Finite-dimensional pointed Hopf algebras over $\PSL_2(q)$}
\label{subsec:psl2}
Frey\-re, Gra\~na and Vendramin studied in \cite{FGV2} 
Nichols algebras associated to conjugacy classes in $\PSL_2(q)$. By the use abelian techniques, they 
proved that 
$\dim \toba(\Oc, \rho) = \infty$ for any $(\rho,V) \in \Irr C_{\PSL_2(q)}(x)$ 
with $x\in \Oc$ and $q$ even, \cite[Proposition 3.1]{FGV1}. For $q$ odd, 
a list of the open cases is given,  see \cite[Theorem 1.6]{FGV2}. 
In the next theorem we use the criterium of type C to discard also
the case when $q\equiv 1 \mod 4$. 

First we recall a lemma that yields that there is no finite-dimensional Nichols algebra
over the conjugacy class of involutions in $\PSL_{2}(7)$. Note that by Proposition \ref{prop:sl2-ss}
this class is kthulhu.

\begin{lema}\cite[Prop. 4.3]{FGV2}\label{lem:sl2-YD-modules}
Let $\Oc$ be the conjugacy class of involutions in $\PSL_{2}(7)$ and $x\in \Oc$. Then, 
$\dim\toba(\Oc,\rho)=\infty$, for every
$\rho\in \Irr C_{\Gb}(x)$.
\end{lema}

\begin{theorem}\label{thm:sl2}
Let $\Oc$ be a semisimple conjugacy class in $\PSL_2(q)$. 
If $\Oc$ is not a semisimple irreducible conjugacy class 
represented in $\SL_{2}(q)$ by $\x = \left(\begin{smallmatrix}
a& \zeta b\\
b & a
\end{smallmatrix}\right)$ with $ab\neq 0$, $\zeta \in \F_{q}^{\times} - \F_{q}^{2}$ and  
$q\equiv 3 \mod 4$, 
then 
$\dim\toba(\Oc,\rho)=\infty$, for every
$\rho\in \Irr C_{\Gb}(x)$.
\end{theorem}

\pf If $q$ is even, then 
$\dim\toba(\Oc,\rho)=\infty$, for every
$\rho\in \Irr C_{\Gb}(x)$ by \cite[Proposition 3.1]{FGV1}, since in this case
$\PSL_{2}(q)=\SL_{2}(q)$. Assume $q$ is odd. 
If $q\equiv 1 \mod 4$, the open cases in \cite[Theorem 1.6]{FGV2} were given
by split semisimple classes which collapse by Proposition \ref{prop:sl2-ss}. 
\epf

\begin{obs}
For $q>3$ and $q\equiv 3 \mod 4$ there are  $\frac{q-3}{4}$ semisimple
irreducible conjugacy classes in $\PSL_{2}(q)$ of size $q(q-1)$, see \cite[Table 3]{FGV2}.
\end{obs}

\section{Finite-dimensional pointed Hopf algebras over $\PSL_n(q)$}\label{sec:psln}
In this last section we show that for infinitely many pairs $(n,q)$, the groups $\SL_{n}(q)$ and 
$\PSL_{n}(q)$ collapse, see Theorem \ref{thm:psl-collapse}. We use both the criteria 
introduced before as well as abelian techniques. We begin by studying irreducible semisimple
conjugacy classes in $\PSL_{3}(q)$.

\subsection{Semisimple classes in  $\PSL_{3}(q)$}\label{sec:psl3-ss}

In this subsection $\Gb=\PSL_3(q)$.  
Let $x$ be an
irreducible semisimple element in $\Gb$. We show 
in Propositions \ref{prop:psl3-odd-ss-kthulhu} and \ref{prop:psl3-even-ss-kthulhu} that
$\Oc_x^{\Gb}$ 
is austere, hence kthulhu. 

If $x$ is semisimple and irreducible, then for every eigenvalue $\eta$ of $\x$ we have
$\eta^{q^2}\neq\eta\neq\eta^q$, $\eta^{(3)_q}=1$, and  
$\x^{(3)_q}=1$, i.e., $\x$ lies in a maximal torus of order ${(3)_q}$. 
Since $Z(\SL_3(q))$ lies in all such tori, $x$ lies in 
the  maximal torus of $\Gb$ of order $\frac{(3)_q}{(3,q-1)}$.

If $x$ is semisimple and not irreducible, then $\x$ lies in a maximal torus whose exponent 
is either $q^2-1$ or $q-1$, whence  $x$ lies in a 
maximal torus of $\Gb$ of order 
$\frac{q^2-1}{(3,q-1)}$ or $\frac{q-1}{(3,q-1)}$.

We need first the following technical lemma. 
\begin{lema}\label{lem:arithmetic}Let $\ell>1$ be odd. If $q\equiv 1 \mod \ell$, then 
$\left(\ell^2, (\ell)_q\right)=\ell$.
\end{lema}
\pf Let $q=1+a\ell$. We show that $(\ell)_q\equiv\ell \mod \ell^2$.
$$
\begin{array}{rl}
(\ell)_q&=\sum_{j=0}^{\ell-1}(1+a\ell)^j=\sum_{j=0}^{\ell-1}
\sum_{i=0}^j\binom{j}{i}a^i\ell^i=\sum_{i=0}^{\ell-1}a^i\ell^i\sum_{j=i}^{\ell-1}\binom{j}{i}\\
\\
&\equiv \sum_{j=0}^{\ell-1}\binom{j}{0}+a\ell\sum_{j=1}^{\ell-1}\binom{j}{1}\mod \ell^2\\
\\
&\equiv \ell +a\ell \binom{\ell}{2} \mod \ell^2 \equiv \ell \mod \ell^2.
\end{array}
$$ Hence, $((\ell)_q,\ell^2)=(\ell,\ell^2)=\ell$. 
\epf

\begin{obs}\label{obs:not9}Let $x=\pi(\x)\in \Gb$, for some 
semisimple irreducible $\x\in\SL_3(q)$. 

\begin{enumerate}
[leftmargin=*,label=\rm{(\alph*)}]
\item \label{item:order} We claim that $(\ord x,6)=1$ and  $\ord x\neq 5$. 
Indeed, $\ord \x$ divides $(3)_q$ which is always odd. 
In addition,  $(3)_q$  is divisible by $3$ only if  
$q\equiv1 \mod 3$ and by Lemma \ref{lem:arithmetic}, it is never divisible by $9$.
Looking at
all the possible values of $q$ modulo $5$ it is easily verified that 
$5\not|(3)_q$. 
\item 
\label{item:index3}
Let $H\leq K\leq \Gb$, with $[K:H]\leq 3$. 
Then, if $x\in K$, we have $x\in H$. Indeed
left multiplication by $x$ induces a permutation of the coclasses of 
$H$ in $K$, which has order $\leq3$. By  
\ref{item:order},
it must be the identity, i.e., $xH=H$.
\item\label{item:tori} If $x^k$ is not irreducible, then $x^k=e$. 
In fact $x^k$ is semisimple and therefore it lies 
in a maximal torus of $\Gb$. The statement follows because 
$\left((3)_q,{(q^2-1)}\right)=(q-1,3)$ and 
$\left(\frac{(3)_q}{(3,q-1)},\frac{(q^2-1)}{(3,q-1)}\right)=1$. 
\end{enumerate}
\end{obs}

\begin{lema}\label{lem:degree3}Let $x=\pi(\x)\in \Gb$ be semisimple 
and irreducible. Assume $q=t^{3l}$. Then
$\Oc_x^{\Gb}\cap \PSL_3(t)=\emptyset$.
\end{lema}
\pf Let $\y\in \SL_3(t)$. If its characteristic polynomial is not irreducible over $\F_t$, 
then  $\pi(\y)\not\in \Oc_x^{\Gb}$. If it is irreducible, 
then its roots are in $\F_{t^3}\subset\F_q$, so 
$\pi(\y)\not\in\Oc_x^{\Gb}$. 
\epf

We prove now that  $\Oc_x^{\Gb}$ 
is austere when $q$ is odd.

\begin{prop}\label{prop:psl3-odd-ss-kthulhu}
Assume that $q$ is odd. Let $K\leq \Gb$ and $x$ be an
irreducible semisimple element in $\Gb$. 
Then one of the following holds:
\begin{enumerate}[leftmargin=*,label=\rm{(\alph*)}]
 \item $K\cap \Oc_x^{\Gb}=\emptyset$;
 \item $K\cap \Oc_x^{\Gb}$ 
is abelian;
 \item there is $y \in K\cap \Oc_x^{\Gb}$ 
such that $K\cap \Oc_x^{\Gb}=\Oc_y^K$.
\end{enumerate}
In particular $\Oc_x^{\Gb}$ 
is austere, hence kthulhu. 
\end{prop}

\pf We proceed by inspection of the different subgroups of 
$\Gb$ as listed in Theorem \ref{th:subgps-psl3-odd}.

If $K$ is as in case
\ref{bloom:point-psl3}, then $K=\PSL_{3}(t)$ for some $t=p^{a}$. 
In this case, we have that if $y\in K\cap \Oc_x^{\Gb}$, then  $K\cap \Oc_x^{\Gb}=\Oc_y^K$.
Indeed, if $(q-1,3)=(t-1,3)$, we apply Lemma \ref{lem:one-gcd-in-pslt}. If $(q-1,3)\neq(t-1,3)$, then
$(q-1,3)=3$ and $(t-1,3)=1$, so  $(3,(3)_t)=1$ and Lemma \ref{lem:irred-one-in-pslt} applies.

If $K$ is as in case \ref{bloom:point-pu3}, then $K=\PSU_{3}(t)$ for some $t=p^{a}$.
Again, we have that if $y\in K\cap \Oc_x^{\Gb}$, 
then  $K\cap \Oc_x^{\Gb}=\Oc_y^K$. For,
if $(q-1,3)=(t+1,3)$, then we apply Lemma \ref{lem:one-in-psu}. If on the other hand, $(q-1,3)\neq(t+1,3)$, then
$(q-1,3)=3$ and $(t+1,3)=1$, $p\neq3$ and $t \equiv 1 \mod 3$. 
Without loss of generality we assume $x=\pi(\x)$ for some $\x\in \SU_3(t)$. 
Let $y=\pi(\y) \in \Oc_x^{\Gb}\cap K$, with $\y,\in \SU_3(t)$. 
Then, for some $\zeta\in\F_q$ with $\zeta^3=1$ 
there holds  $\y\in \Oc_{\zeta\x}^{\SL_3(q)}\cap \SU_3(t)$.
There are three conjugacy classes of  maximal tori in $\SU_3(t)$, 
with exponent $t+1$, $t^2-1$, which both divide $q-1$, 
and $t^2-t+1$. Since $y$ and $x$ are irreducible, 
$\ord(\zeta \x)=\ord(\y)$ and $\ord(\x)$ divide $t^2-t+1$. 
But then, $\zeta=\zeta^{t^2-t+1}=1$. Hence, 
$\y\in \Oc_{\x}^{\SL_3(q)}\cap \SU_3(t)=\Oc_{\x}^{\SU_3(t)}$ 
where equality follows because 
the centralizer of a semisimple element in $\SL_3(\kk)$ is connected. Thus, $y\in \Oc_x^K$.

If $K$ is as in case \ref{bloom:point-psl3-extension},
then there is a subgroup containing  $\PSL_{3}(t)$ with index 3 and
$t\equiv 1 \mod 3$. In this case, 
$K\cap \Oc_x^{\Gb}=\emptyset$. Indeed, if 
$y\in \Oc_x^{\Gb}$, then by Lemma \ref{lem:degree3}, $y\not\in \PSL_3(t)$, and 
consequently
$y\not\in K$ by Remark \ref{obs:not9} \ref{item:index3}.

If $K$ is as in case
\ref{bloom:point-pu3-extension},   
then there is a subgroup containing  $\PSU_{3}(t)$ with index 3 and
$t\equiv 2 \mod 3$. As before, it follows that 
$K\cap \Oc_x^{\Gb}=\emptyset$. 
If $y\in \Oc_x^{\Gb}$, then by Lemma \ref{lem:degree3}, 
$y\not\in \PSU_3(t)\subset \PSL_3(t^2)$. Thus, by Remark \ref{obs:not9} 
\ref{item:index3}, $y\not\in K$.

If $K$ is as in case
\ref{bloom:point-psl2}, then $K= \PSL_{2}(t)$ or $\PGL_{2}(t)$ for some $t=p^{a}\neq 3$.
In this case, $K\cap \Oc_x^{\Gb}=\emptyset$. Indeed, 
the order of $K$ divides $t(t^2-1)$ which in turn divides $t(q^2-1)$. If $x\in K$ were
irreducible, then $\ord(x)$ would be a divisor of
$(t(q^2-1),(3)_q)=(3,q-1)$. Hence, Remark \ref{obs:not9} \ref{item:order} applies.

If $K$ is as in cases
\ref{bloom:point-psl2(5)}, \ref{bloom:point-psl2(7)}, \ref{bloom:point-a67}, 
\ref{bloom:point-a6}, then the order of any element in $K$ lies in $\{2,3,4,5,7\}$. By Remark \ref{obs:not9} \ref{item:order}, 
if $y\in K\cap \Oc_x^{\Gb}$, then $\ord y=7$ and $K$ 
is either ${\mathbb A}_7$ or $\PSL_2(7)$. 
In addition, $7$ divides $(3)_q$ only if 
$q\equiv 2 \mod 7$ or $q\equiv 4 \mod 7$. For both choices of $K$ there are exactly 
two classes of elements of order $7$, one containing 
$y, y^2=y^q, y^4=y^{q^2}$, and the other containing $y^3, y^5, y^6$. 
We  show that $\Oc_y^K=\Oc_y^{\Gb}\cap K$. Assume 
this is not the case. Then, $y^{3}\in \Oc_y^{\Gb}$, so $x^3\in \Oc_x^{\Gb}$.
If $x=\pi(\x)$, then the only powers of $\x$ lying in 
$\Oc_{\x}^{\SL_3(q)}$ are $\x,\x^q$ and $\x^{q^2}$.
By looking at the order of $\x$ (which can be $7$ or $21$), 
and of $\x^3$ (which is always $7$)
we necessarily have $\x^3\in \Oc_{z\x}^{\SL_3(q)}$
for some $1\neq z\in Z(\SL_3(q))$, a third root of $1$, with $\x^7=z^{-1}$, 
and $q\equiv 1 \mod 3$. 
This is impossible. Indeed, let $\zeta, \zeta^q$ and $\zeta^{q^2}$ 
be the eigenvalues of $\x$. They are primitive 
$21$-th roots of $1$. Then, the eigenvalues of $\x^3$ 
are $\zeta^3, \zeta^{3q}$ and $\zeta^{3q^2}$, whereas the 
eigenvalues of $z\x$ are $z\zeta, z\zeta^q$ and $z\zeta^{q^2}$. 
A direct verification using that $q\equiv 2\mod 7$ or 
$q\equiv 4 \mod 7$
shows that $\zeta^3$ cannot be equal to any of $z\zeta, z\zeta^q$ and $z\zeta^{q^2}$. 
   
If $K$ is as in case \ref{bloom:point-cyclic-normal}, then
$K$ contains a cyclic $p$-regular normal subgroup $H=\langle h\rangle$ of index $\leq 3$, which
is contained in a maximal torus $S$. 
Let $y\in K$. If $y$ is irreducible in $\Gb$,
 then by Remark \ref{obs:not9} \ref{item:order}, 
\ref{item:index3}, $y\in H$ and $H$ is contained in the torus $S$ of order  
$\frac{(3)_q}{(3, (3)_q)}$. Therefore, 
 $\Oc_x^{\Gb}\cap K=\Oc_x^{\Gb}\cap H$ is abelian. 

 If $K$ is as in case 
\ref{bloom:point-diagonal-normal}, then $K$ contains a diagonal normal 
subgroup $L$ such that $K/L $ is isomorphic to a subgroup of $ \st$. In this case, we also have that
$\Oc_x^{\Gb}\cap K=\emptyset$. Assume on the contrary that $x\in K$. 
The order of its coclass $xL$ in 
$K/L$ is $1,2$, or $3$. It cannot be $1$ because $x$ is 
irreducible hence not diagonal, and it cannot be $2$ nor $3$ 
by Remark \ref{obs:not9} \ref{item:order}. 

If $K$ is as in case 
\ref{bloom:point-diagonal-diese}, then $K$ contains a normal abelian $p$-subgroup $H$ such 
that $K/H$ is isomorphic to a subgroup of $\GL_{2}(q)$. In particular, 
$|K|$ divides $p^N|\GL_2(q)|=p^Nq(q-1)(q^2-1)$, for some $N>0$. 
As
$\left(p^N(q-1)(q^2-1), \frac{(3)_q}{(3,q-1)}\right)=1$, it follows that
$K\cap \Oc_x^{\Gb}=\emptyset$. 

If $K$ is as in case 
\ref{bloom:point-diagonal-abelian33}, then 
$K$ contains a normal abelian subgroup $H$ of type $(3,3)$ such that  
$K/H$  is isomorphic to a subgroup of $\SL_{2}(3)$.
We  show that $\Oc_x^{\Gb}\cap K=\emptyset$. 
Assume that $y\in K$ is semisimple and irreducible. 
We look at the coclass $yH$ in $K/H$. 
Then $\ord yH$ is either $1,2,3,$ or $4$, whence
by Remark \ref{obs:not9} \ref{item:order} we have 
$y\in H$, a $3$-group, which is impossible.  

Finally, if $K$ is as in case
\ref{bloom:point-diagonal-quaternion}, then 
$K$ contains a normal abelian subgroup $H$ of type $(3,3)$ such that 
$K/H$  is isomorphic to a subgroup of the quaternion group $\Q$ of order $8$.
Hence, 
$\Oc_x^{\Gb}\cap K=\emptyset$ 
by Remark \ref{obs:not9} \ref{item:order}.
\epf

We prove now that  $\Oc_x^{\Gb}$ 
is austere when $q$ is even.

\begin{prop}\label{prop:psl3-even-ss-kthulhu}
Assume that $q=2^m$. 
Let $x$ be an irreducible 
semisimple element in $\Gb$. Then, for $y\in \Oc_x^{\Gb}$ we have either $xy=yx$ or 
$\Oc_x^{\langle x,y\rangle}=\Oc_y^{\langle x,y\rangle}$.
In particular $\Oc_x^{\Gb}$ is austere, hence kthulhu. 
\end{prop}
\pf Let $x,y\in  \Oc_x^{\Gb}$, with $xy\neq yx$. If $K:=\langle x,y\rangle\neq \Gb$, 
then $K$ lies in a proper maximal subgroup $M_1$ of $\Gb$.  
We analyse the different possibilities for $M_1$ listed in  Theorem \ref{th:subgps-psl3-even}, 
from which we adopt notation.  

Since $x$ is a semisimple irreducible element, 
$\ord x$ is coprime with the order of the 
groups in  \ref{hartley:normalizer-tori}, thus this case is not possible. Further, 
$\ord x$ is different from the 
order of any element in  \ref{hartley-a6}
by Remark \ref{obs:not9} \ref{item:order}. 
If $M_1$ were 
as in \ref{hartley:normalizer-maximal torus}, then by 
Remark \ref{obs:not9} \ref{item:index3}, $x$ and $y$ would lie 
in the same maximal torus, and they would commute.  
Hence, the possible 
cases are \ref{hartley-psl}, 
\ref{hartley-pgl}, \ref{hartley-psu} and \ref{hartley-pgu}, 
that is $M_{1}$ is isomorphic to $\PSL_3(t)$ or $\PSU_{3}(t)$ for some $t$, or contains one of them with index $3$.
We analyze these cases further.

\smallbreak
\noindent{\emph Claim 1.} $M_1$ is  either 
$\PSL_3(t)$, where $q$ is a prime power of $t$ and $q\neq t^3$, or $\PSU_3(t)$. 
The latter may occur only if $m$ is even.

If $x,y\in \PGL_3(t)$ with $q = t^3$, then by  
Remark \ref{obs:not9} \ref{item:index3}, $x, y\in \PSL_3(t)$,
which is impossible by Lemma \ref{lem:degree3}, so case 
\ref{hartley-pgl} is excluded. Similarly, if $x,y\in \PGU_3(t)$, with $q = t^6$,  then 
$x,y\in \PSU_3(t)\leq\PSL_3(t^2)$, impossible by Lemma \ref{lem:degree3}. Thus,
case  \ref{hartley-pgu} is also excluded.

Now we proceed inductively looking at the maximal subgroups of $M_{1}$ as above.

\smallbreak
\noindent{\emph Claim 2.} $K$ 
is either  $\PSL_3(t')$, for some $t'=2^b$ and $b|m$, $3b\not|\ m$, 
or $\PSU_3(t')$, for $t'=2^{c}$ and $c|2m$, 
$3c\not|\ 2m$.

If $K=M_1$ the claim is trivial. Otherwise, $K \le M_2$ where $M_2$ 
is a maximal subgroup of $M_1$. If $M_1=\PSL_3(t)$ we argue as in Case 1. 
If $M_1=\PSU_3(t)$,
then we claim that $M_2=\PSU_3(t')$ for $t$ an odd prime power of $t'$.
We analyse the different possibilities for $M_2$ listed in  
Theorem \ref{th:subgps-psu3-even},
from which we adopt notation. The groups as in 
\ref{hartley:psu-normalizer-tori} or \ref{hartley-36}
are discarded because their order is coprime with $\ord x$. 
The groups as in \ref{hartley:psu-normalizer-maximal torus}
 may not occur by Remark \ref{obs:not9} \ref{item:index3} 
and the  noncommutativity of $x$ and $y$.
If $x$ would lie in a subgroup as in \ref{hartley-pgu-psu}, 
by Lemma \ref{lem:degree3} it would lie in 
$\PSU_3(2^l)\leq \PSL_3(2^{2l})$ with $t=2^a$, $6l=2a|m$. 
By Lemma \ref{lem:degree3} we have a contradiction.
Thus, the only possible case is \ref{hartley-36} and we have the claim.

\smallbreak
\noindent{\emph Claim 3.} $\Oc_x^K=\Oc_x^{\Gb}\cap K=\Oc_y^K$.

Let us observe that $x,y$ are again irreducible in $K$. 
If $K=\PSL_3(t)$ for some $t$ the argument in the proof of 
Proposition \ref{prop:psl3-odd-ss-kthulhu} case \ref{bloom:point-psl3} gives the claim. 
If $K=\PSU_3(t)$, we  argue as in the proof of 
Proposition \ref{prop:psl3-odd-ss-kthulhu}  case \ref{bloom:point-pu3}.
\epf

\subsection{Finite-dimensional pointed Hopf algebras over $\PSL_n(q)$}\label{subsec:psln}
In this last subsection we assume that $\Gb=\PSL_n(q)$, $n>2$. 
We show that for infinitely many pairs $(n,q)$, the groups $\SL_{n}(q)$ and 
$\PSL_{n}(q)$ collapse.

Let $n=2^ab$, $q=1+2^c d$ 
with $(b,2)=(d,2)=1$ and let $\mathcal{G}_{ss}$ be the set of pairs $(n,q)$ 
with $n\in \N$, $n>2$ and $q=p^{m}$, $p$ a prime such that 
one of the following hold:
\begin{enumerate}[leftmargin=*,label=\rm{(\alph*)}]
\item $n$ is odd;
\item $q$ is even;
\item $0<a<c$;
\item $a=c>1$.
\end{enumerate}

A direct computation similar to the one in the proof of Lemma \ref{lem:arithmetic} shows that if $(n,q) \in \mathcal{G}_{ss}$, then 
$\frac{(n)_q}{(q-1,n)}$ is odd.

\begin{lema}\label{lem:irred-psll-typeB}Let $(n,q) \in \mathcal{G}_{ss}$ and let
 $x\in \PSL_n(q)$ be an irreducible semisimple element. 
 Then, for every $\rho\in {\rm Irr}{C_{\Gb}(x)}$ we have 
$\dim \toba(\Oc_x^{\Gb},\rho)=\infty$. The same statement holds 
for $G=\SL_n(q)$.
\end{lema}

\pf  By Lemma \ref{lem:three-powers}, 
$x\neq x^q$, $x\neq x^{q^2}$ and $x^q, x^{q^2}\in \Oc_x^{\Gb}$. 
In addition the order of the  maximal torus 
containing $x$ is  $\frac{(n)_q}{(q-1,n)}$ which is odd 
under our assumptions. Then, Lemma \ref{lem:ab-techniques-real} (c) applies. 
\epf

\begin{prop}Let $\Oc$ be a conjugacy class in $\Gb=\PSL_3(2)$. If $\Oc$ 
is not unipotent of type $(3)$, then
$\dim \toba(\Oc_x^{\Gb},\rho)=\infty$
for every $\rho\in {\rm Irr}\,{C_{\Gb}(x)}$.
\end{prop}
\pf This follows from Propositions  \ref{prop:sln-nonss}, \ref{prop:reduction}  
and Lemma \ref{lem:irred-psll-typeB}.
\epf

We end the paper with the following result. Recall the definition of
the set $\mathcal{F}$ of pairs $(n,q)$ given in the Introduction.

\begin{theorem}\label{thm:psl-collapse} If $(n,q)\in \mathcal{F}$,
then $\PSL_n(q)$  and $\SL_n(q)$ collapse. 
\end{theorem}

\pf Let $\Oc$ be a conjugacy class in $\Gb$. If $\Oc$ is not semisimple, 
then we apply Proposition \ref{prop:sln-nonss}. If $\Oc$ is semisimple,
then we use Proposition \ref{prop:reduction} and Lemma \ref{lem:irred-psll-typeB}.
\epf

\begin{obs}\label{obs:abelian}Assume $n=2^ab$ and $q=1+2^c d$ 
with $(b,2)=(d,2)=1$ and $a>c>0$ or else $a=c=1$. 
Let $x$ be an irreducible semisimple element in $\Gb$. If $\ord x$ is odd, 
then $\dim\toba(\Oc_x,\rho)=\infty$ for all $\rho\in\Irr C_{\Gb}(x)$ 
by Lemmata \ref{lem:three-powers} and \ref{lem:ab-techniques-real} (c). Hence, 
the only potentially non-collapsing classes in $\Gb$ are the classes of 
semisimple irreducible elements of even order. In this case such elements always exist: 
for instance,  any generator of a  maximal torus $S$ of order 
$\frac{(n)_q}{(q-1,n)}$, which is even under our assumptions. 
Any semisimple irreducible element is conjugate to an element in $S$.

Let $x=\pi(\x)$ be an irreducible semisimple element in $S$. Its centralizer in $\Gb$ is
$$C_{\Gb}(x)=\{\pi(\y)\in\Gb~|~\y\x\y^{-1}=\lambda\x,\,\lambda\in\F_q\}.$$
If for some $\y\in G$ there holds $\y\x\y^{-1}=\lambda \x$, 
then $\lambda\x\in\Oc\cap C_{G}(\x)=\{\x^{q^i},\,i\in\I_{0,n-1}\}$. 
By Lemma \ref{lem:three-powers} this implies $\lambda=1$ and consequently 
$\y\in C_G(\x)$. Thus, $C_{\Gb}(x)=\pi(C_G(\x))=S$ is cyclic, so all its irreducible 
representations are $1$-dimensional. 
By Lemma \ref{lem:ab-techniques-real} (c), if $\dim\toba(\Oc_x,\rho)<\infty$, 
then $\rho(x)=-1$. In this case the study of abelian 
subracks of $\Oc$ is not effective for determining the dimension of the associated 
Nichols algebras. 
\end{obs}

\subsection*{Acknowledgements} We thank Robert Guralnick 
for conversations on subgroups generated by two elements in finite simple groups  and 
Jay Taylor for interesting discussions on semisimple classes.


\begin{thebibliography}{AFGV2}

\bibitem[ACG1]{ACG-I}  N. Andruskiewitsch, G. Carnovale, G. A.
Garc\'ia.
\emph{Finite-dimensional pointed Hopf algebras over finite simple groups of Lie type  I. 
Non-semisimple classes 
in $\PSL_n(q)$}, 
J. Algebra, to appear,  doi:10.1016/j.jalgebra.2014.06.019

\bibitem[ACG2]{ACG-II}  \bysame
\emph{Finite-dimensional pointed Hopf algebras over finite simple 
groups of Lie type  II. Unipotent classes in symplectic groups}, 
\texttt{arXiv:1412.7397}. Commun. Contemp. Math., to appear,  DOI: 10.1142/S0219199715500534. 


\bibitem[AF]{AF2} N. Andruskiewitsch and F. Fantino,
\emph{On pointed Hopf algebras associated with alternating and
dihedral groups}, Rev. Uni\'on Mat. Argent. 48-3, (2007), 57-71.


\bibitem[AFGaV1]{AFGV-simple}  N. Andruskiewitsch, F. Fantino, G. A.
Garc\'ia and  L. Vendramin,
\emph{On Nichols algebras associated to simple racks},
Contemp. Math. \textbf{537} (2011), 31--56.

\bibitem[AFGV1]{AFGV-ampa} N. Andruskiewitsch, F. Fantino, M.
Gra\~na and  L. Vendramin, \emph{Finite-di\-mensional pointed
Hopf algebras with alternating groups are trivial},
Ann. Mat. Pura Appl. (4),  \textbf{190}  (2011), 225--245.

\bibitem[AFGV2]{AFGV-espo} \bysame,
\emph{Pointed Hopf algebras over the sporadic simple groups}.
J. Algebra \textbf{325} (2011), pp. 305--320.


\bibitem[AG]{AG-adv} N. Andruskiewitsch and M. Gra\~na,
    \emph{From racks to pointed Hopf algebras},
Adv. Math.  \textbf{178}  (2003), 177--243.

\bibitem[AHS]{AHS} N. Andruskiewitsch, I. Heckenberger and H.-J. Schneider,
\emph{The Nichols algebra of a semisimple Yetter-Drinfeld module},
Amer. J. Math.  \textbf{132}, no. 6, 1493--1547.

\bibitem[AZ]{AZ}  N. Andruskiewitsch and  S. Zhang,
\emph{On pointed Hopf algebras associated to some conjugacy
classes in $\mathbb S_n$}. Proc. Amer. Math. Soc. \textbf{135}  (2007),
2723-2731.

\bibitem[B]{B} D. M. Bloom,
\emph{The Subgroups of PSL(3, q) for odd q}.
Trans. Am. Math. Soc. \textbf{127} (1967), 150--178.


\bibitem[Bou]{Bou} {N. Bourbaki},
{\it Groupes et alg\`{e}bres de Lie Chap. IV,V,VI}, Hermann, Paris, 1968.

\bibitem[BGK]{BGK}
T. Breuer; R. Guralnick; W. M. Kantor. \emph{Probabilistic generation of finite simple groups. II}. 
J. Algebra \emph{320} (2008), 443--494.







\bibitem[CH]{CH} M. Cuntz and I. Heckenberger,
\emph{Finite Weyl groupoids}. J. reine angew. Math., to appear.


\bibitem[D]{D} J. Dieudonne, 
\emph{La geometrie des groupes classiques}. 
Ergebnisse der Mathematik und ihrer Grenzgebiete. 
Neue Folge ; vol. 5. Berlin, Springer, 1963.

\bibitem[FaV]{FaV} {F. Fantino} and {L. Vendramin},
\textit{On twisted conjugacy classes of type D in sporadic simple groups}.
 Contemp. Math. \textbf{585} (2013) 247-259.

\bibitem[FGV1]{FGV1} S. Freyre, M. Gra\~na and L. Vendramin,
\emph{On Nichols algebras over $\SL(2,q)$ and $\GL(2,q)$}. J. Math. Phys. \textbf{48}, (2007) 123513.


\bibitem[FGV2]{FGV2} \bysame,  
\textit{On Nichols algebras over $\PGL (2,q)$ and $\PSL(2,q)$}.
J. Algebra Appl. \textbf{9} (2010), 195--208.


%
%

\bibitem[GK]{GK}
R. Guralnick; W. M. Kantor. \emph{Probabilistic generation of finite simple groups}. 
 Special issue in honor of Helmut Wielandt. J. Algebra 234 (2000),  743--792.


 \bibitem[GM]{GM}
R. Guralnick; G. Malle.  \emph{Products of conjugacy classes and fixed point spaces}. 
J. Amer. Math. Soc. 25 (2012), 77-121.
 
 \bibitem[Ha]{hartley}
R. W. Hartley. \emph{Determination of Ternary Collineation Groups Whose Coefficients Lie in the $GF(2n)$}
Ann. Math. \textbf{27}(2) (1925), 140--158. 
 
\bibitem[H]{H-all}
I. Heckenberger, \textit{Classification of arithmetic root systems},
Adv. Math. {\bf 220} (2009) 59--124.

\bibitem[HS]{HS} I. Heckenberger, and  H.-J. Schneider,
\textit{Nichols algebras over groups with finite root system of rank two {I}}.
J. Algebra \textbf{324}(11), 3090--3114 (2010).

\bibitem[HV]{HV}
I. Heckenberger and L. Vendramin, 
\textit{The classification of Nichols algebras with finite root system of rank two},
J. Europ. Math. Soc., to appear.

\bibitem[Hu]{Hu} J. E. Humphreys,
\emph{Conjugacy classes in semisimple algebraic groups},
Amer. Math. Soc., Providence, RI, 1995.




\bibitem[LS]{LS}M. W. Liebeck, G. M. Seitz,
\emph{Unipotent and Nilpotent Classes in Simple Algebraic
Groups and Lie Algebras},
Amer. Math. Soc. Providence, RI, (2012).

\bibitem[MaT]{MT} G. Malle and D. Testerman,
\emph{Linear Algebraic Groups and Finite Groups of Lie Type},
Cambridge Studies in Advanced Mathematics \textbf{133} (2011).


%
%
%

\bibitem[Su]{suzuki}M. Suzuki,
\emph{Group Theory I},
Grundlehren der mathematischen Wissenschaften 247; Springer (1982).



\end{thebibliography}
\end{document}